\documentclass[11pt]{article}
\usepackage[margin=1in]{geometry}
\usepackage{amsmath,amssymb,amsthm,mathtools}
\usepackage{enumitem}
\usepackage{graphicx}
\usepackage{todonotes}
\usepackage{tikz}
\usetikzlibrary{arrows.meta}
\usepackage{booktabs}
\usepackage[square,numbers]{natbib}
\usepackage[colorlinks=true,linkcolor=blue,citecolor=blue,urlcolor=blue]{hyperref}
\newcommand*{\doi}[1]{doi:~\href{https://doi.org/#1}{\nolinkurl{#1}}}

\newtheorem{theorem}{Theorem}[section]
\newtheorem{proposition}[theorem]{Proposition}
\newtheorem{lemma}[theorem]{Lemma}

\theoremstyle{definition}

\theoremstyle{remark}
\newtheorem{remark}[theorem]{Remark}

\newcommand{\X}{\mathcal X}
\newcommand{\Y}{\mathcal Y}
\newcommand{\PP}{\mathcal P}
\newcommand{\E}{\mathbb E}
\newcommand{\R}{\mathbb R}

\newcommand{\dX}{d_{\X}}
\newcommand{\dY}{d_{\Y}}
\newcommand{\muh}{\hat\mu_n}
\newcommand{\eps}{\varepsilon}

\DeclareMathOperator{\diam}{diam}
\DeclareMathOperator{\Hess}{Hess}

\title{The Rate of Convergence of Empirical $(p,q)$-Gromov--Wasserstein Distances}
\author{%
  Pui Kuen Leung\thanks{Department of Statistics, University of Chicago. Email: \texttt{pkl@uchicago.edu}.}
  \and
  Riku Okada\thanks{Committee on Computational and Applied Mathematics, University of Chicago. Email: \texttt{rikuo@uchicago.edu}.}
  \and
  Samuel Lok-Hei Wong\thanks{Committee on Computational and Applied Mathematics, University of Chicago. Email: \texttt{slhwong@uchicago.edu}.}
}
\date{\today}
\makeatletter
\g@addto@macro\@maketitle{%
  \begingroup
  
  \footnotetext[0]{The authors are listed in alphabetical order.}
  \endgroup
}
\makeatother

\begin{document}
\maketitle
\begin{abstract}
We study empirical convergence rates for $(p,q)$-Gromov--Wasserstein (GW) distances $D_{p,q}$ \citep{sturm2023space} for arbitrary real exponents $p>0$ and $q\ge1$. Let $\mu$ and $\nu$ be probability measures on compact metric spaces $\X$ and $\Y$, and let $\hat\mu_n$ and $\hat\nu_m$ be their empirical measures, constructed from independent i.i.d.\ samples of sizes $n$ and $m$. If the spaces $\X,\Y$ are compact subsets of the Euclidean spaces $\R^{d_x}$ and $\R^{d_y}$, respectively, we obtain an upper bound of order $(n\wedge m)^{-(p\wedge2)/((d_x\wedge d_y)\vee2(p\wedge2))}$ on the plug-in error of $D_{p,q}^p(\hat\mu_n,\hat\nu_m)$ in expectation, up to a logarithmic factor at the critical dimension $d_x\wedge d_y=2(p\wedge2)$. On Euclidean balls, this rate is sharp up to the logarithmic factor. More generally, we show that the same upper bound holds when one space is a compact subset of a complete, connected $d$-dimensional Riemannian manifold and the other has polynomial covering numbers, with $d_x\wedge d_y$ replaced by $d$.
Earlier works on $(2,2)$-GW, $(2,1)$-GW and other even-order cases rely on variational representations of GW distances that are specific to their exponents \citep{zhang2024gromov,leung2026sample,paliy2026empirical}. We instead reduce the plug-in error to an analysis of empirical processes through local duality, which may be of independent interest.


\end{abstract}

\section{Introduction}\label{sec:intro}

Let $(\X,\dX,\mu)$ and $(\Y,\dY,\nu)$ be compact metric measure spaces. Given $p,q>0$, the
$(p,q)$-Gromov--Wasserstein (GW) distance $D_{p,q}(\mu,\nu)$, also called the $L^{p,q}$-distortion
distance \citep{sturm2023space}, is defined as
\begin{equation}\label{eq:GW}
  D_{p,q}(\mu,\nu)=\Bigl(\inf_{\pi\in\Pi(\mu,\nu)}\iint_{(\X\times\Y)^2}\bigl|\dX(x,x')^q-\dY(y,y')^q\bigr|^p\,d\pi(x,y)\,d\pi(x',y')\Bigr)^{1/p},
\end{equation}
where $\Pi(\mu,\nu)$ is the set of couplings of $\mu$ and $\nu$. The
$(p,q)$-GW distance extends the
$p$-GW distance \citep{memoli2011gromov}, introduced by Mémoli as an $L^p$ relaxation of the
Gromov--Hausdorff distance, by raising the pairwise distances in the integrand to a power
$q$. As is known for the $p$-GW distance \citep[Theorem~5.1(a)]{memoli2011gromov}, $(p,q)$-GW defines a metric on the
collection of isomorphism classes of compact metric measure spaces when $p\ge1$.

Given the empirical measures
\begin{equation}\label{eq:emp}
  \hat\mu_n=\frac1n\sum_{i=1}^n\delta_{X_i},\qquad \hat\nu_m=\frac1m\sum_{j=1}^m\delta_{Y_j}
\end{equation}
constructed from independent i.i.d.\ samples $X_1,\dots,X_n\sim\mu$ and $Y_1,\dots,Y_m\sim\nu$, a
slight modification of the proof of M{\'e}moli~\citep[Theorem~5.1(e)]{memoli2011gromov} shows that
$D_{p,q}(\hat\mu_n,\hat\nu_m)\to D_{p,q}(\mu,\nu)$ almost surely as $n,m\to\infty$, for all
$p,q>0$. The aim of this paper is to quantify the
rate of this convergence. We are interested in nonasymptotic bounds of the form
\begin{equation}\label{eq:form}
  \E\bigl|D_{p,q}^p(\hat\mu_n,\hat\nu_m)-D_{p,q}^p(\mu,\nu)\bigr|\le C\,r(n\wedge m),
\end{equation}
where $C$ depends only on $p$, $q$ and the two spaces, and $r(n\wedge m)\to0$ as $n\wedge m\to\infty$.

When $\X$ and $\Y$ are finite,
Rioux et al.~\citep{rioux2026discrete} recently established a variational
representation for general $(p,q)$-GW and obtained bounds of the
form~\eqref{eq:form} for equal sample sizes $n=m$, with $r(n)=n^{-1/2}$.

For Euclidean compact subsets $\X\subset\R^{d_x}$
and $\Y\subset\R^{d_y}$, with $p,q\ge1$, one can obtain a bound of the
form~\eqref{eq:form} by comparing $(p,q)$-GW with $p$-Wasserstein distances
\citep[Lemma~2.1(i)]{kato2025convergence},
\citep[Lemma~4.4]{zhang2024gromov}
and applying classical rates for empirical measures in Wasserstein
distance~\citep{fournier2015rate}.
However, such an argument only yields bounds of the form~\eqref{eq:form}
with $r(n\wedge m)=(n\wedge m)^{-1/((d_x\vee d_y)\vee2p)}$,
up to a logarithmic factor when $d_x\vee d_y=2p$.
Here the exponent depends on the larger dimension $d_x\vee d_y$.

Sharper bounds, in which the exponent depends on the dimensions only
through their minimum $d_\star=d_x\wedge d_y$, have been obtained for
$(p,q)=(2,1)$ \citep{leung2026sample}, for $(2,2)$
\citep{zhang2024gromov}, and for $(2r,2k)$ with $r,k$
positive integers \citep{paliy2026empirical}.\footnote{Although these works consider equal sample sizes $n = m$,
the bounds extend readily to unequal sample sizes; see for instance
\citep[Section~2.3, eq.~(7)]{kato2025convergence}
for the unequal-sample version when $(p,q)=(2,2)$.}
To obtain this dependence on $d_\star$,
these works use arguments based on the lower complexity adaptation principle
from optimal transport \citep{hundrieser2024empirical}, where the reduction to
optimal transport relies on variational representations of $D^p_{p,q}$
specific to these exponents
(\citep[Corollary~4.1]{zhang2024gromov};
\citep[Theorem~4.1 and Corollary~4.2]{leung2026sample};
\citep[Theorem~2.4]{paliy2026empirical}).

In this paper, we extend these $d_\star$-dependent bounds to all real $p>0$ and $q\ge1$.
We also go beyond the Euclidean setting by proving an upper bound when
one space is a compact subset of a Riemannian manifold and the other
has polynomial covering numbers.
We reduce the problem to optimal transport using local duality,
without a variational representation of $D_{p,q}^p$.
We also establish a lower bound for Euclidean balls, which shows that
the upper bound we obtain is sharp up to a logarithmic factor.

\subsection{Setting and main results}\label{sec:setting}

We now state upper and lower bounds in the Euclidean setting, followed by
a Riemannian extension of the upper bound. Throughout, we write $a\wedge b=\min(a,b)$ and $a\vee b=\max(a,b)$.
For a metric space $S$, let $\PP(S)$ denote the set of Borel probability measures on $S$.

For $s>0$, $d\ge1$ and $n\ge2$, set
\begin{equation}\label{eq:psi}
  \psi_{d,s}(n)=\begin{cases}n^{-1/2},&d<2s,\\ n^{-1/2}\log n,&d=2s,\\ n^{-s/d},&d>2s.\end{cases}
\end{equation}
Thus $\psi_{d,s}(n)$ is of polynomial order $n^{-s/(d\vee 2s)}$, up to a logarithmic factor when $d=2s$.

\medskip\noindent\emph{Euclidean setting.}
We first consider the case where both spaces are compact subsets of Euclidean spaces. The following theorem gives an upper bound for $p>0$, $q\geq1$, where the rate exponent depends on the smaller of the ambient Euclidean dimensions, together with a lower bound on Euclidean balls. The latter applies to \emph{all} $p,q > 0$, including the case $q < 1$ omitted for the upper bound, and shows that the upper bound cannot in general be
improved, apart from the logarithmic factor at the critical dimension.
\begin{theorem}[Euclidean bounds]\label{thm:euclid}
\leavevmode
\begin{enumerate}[label=\textup{(\alph*)},ref=\thetheorem(\alph*)]
\item Let $p > 0$, $q\geq 1$. Let $\X\subset\R^{d_x}$ and $\Y\subset\R^{d_y}$ be compact, with $d_x,d_y\ge1$, and set
$d_\star:=d_x\wedge d_y$. Then, for all $\mu\in\PP(\X)$, $\nu\in\PP(\Y)$ and $n,m\ge2$,
\begin{equation}\label{eq:euclidupper}
  \E\bigl|D_{p,q}^p(\hat\mu_n,\hat\nu_m)-D_{p,q}^p(\mu,\nu)\bigr|
  \le C\bigl[\psi_{d_\star,p\wedge2}(n)+\psi_{d_\star,p\wedge2}(m)\bigr],
\end{equation}
where $C$ depends only on $p$, $q$, $d_x$, $d_y$, $\X$ and $\Y$. In particular, the bound is of order $(n\wedge m)^{-(p\wedge2)/(d_\star\vee 2(p\wedge2))}$, with an additional factor $\log(n\wedge m)$ when $d_\star=2(p\wedge2)$.
\item\label{thm:lower} Let $p,q>0$ and $d_x,d_y\ge1$, and set $d_\star:=d_x\wedge d_y$.
For the closed unit balls $B_1^{d_x}(0)$ and $B_1^{d_y}(0)$, there are constants $c>0$ and $n_0$,
depending only on $p$, $q$ and $d_\star$, such that, for all $n,m\ge n_0$,
\begin{equation}\label{eq:lowerintro}
  \sup_{\substack{\mu\in\PP(B_1^{d_x}(0))\\\nu\in\PP(B_1^{d_y}(0))}}\E\bigl|D_{p,q}^p(\hat\mu_n,\hat\nu_m)-D_{p,q}^p(\mu,\nu)\bigr|
  \ge c\bigl(n^{-(p\wedge2)/(d_\star\vee2(p\wedge2))}+m^{-(p\wedge2)/(d_\star\vee2(p\wedge2))}\bigr).
\end{equation}
\end{enumerate}
\end{theorem}
When $p\ge2$, the upper bound is of order $(n\wedge m)^{-2/(d_\star\vee4)}$, up to a logarithmic factor when $d_\star=4$. This matches the rates previously obtained for $(p,q)=(2,1)$, $(2,2)$, and $(2r,2k)$ with integers $r,k\ge1$ \citep{leung2026sample,zhang2024gromov,paliy2026empirical}. When $p<2$, the upper bound~\eqref{eq:euclidupper} is of order
$(n\wedge m)^{-p/(d_\star\vee2p)}$, up to a logarithmic
factor when $d_\star=2p$. This is slower than the rate for $p\ge2$ when $d_\star>2p$. For instance, for the $1$-GW distance $D_{1,1}$, we
obtain the rate $(n\wedge m)^{-1/(d_\star\vee2)}$, up to a logarithmic factor when $d_\star=2$. This has the same order as the classical upper bound for the
convergence of empirical measures in $1$-Wasserstein distance
in dimension $d_\star$~\citep[Theorem~1]{fournier2015rate}.

The rate exponent and critical dimension in~\eqref{eq:euclidupper} are
independent of $q$ (see Step~2 of the proof overview in
Section~\ref{sec:overview}), although the leading constant depends on $q$.

The lower bound in Theorem~\ref{thm:lower} generalizes those obtained in prior work
\citep{zhang2024gromov,leung2026sample}, with one subtle clarification for $(p,q)=(2,2)$. In
Zhang et al.~\citep[Theorem~4.2]{zhang2024gromov}, the lower bound is stated for arbitrary compact sets
$\X\subset\R^{d_x}$ and $\Y\subset\R^{d_y}$, with the supremum taken over $\mu\in\PP(\X)$ and
$\nu\in\PP(\Y)$. However, when $d_\star>4$, their argument involves the construction of certain probability measures on closed balls of nonzero radius, which need not be contained in a given compact set. Our proof uses a similar construction, which motivates stating the lower bound~\eqref{eq:lowerintro} for the closed unit balls $B_1^{d_x}(0)$ and $B_1^{d_y}(0)$.

\medskip\noindent\emph{Riemannian setting.}
Our next result generalizes the upper bound~\eqref{eq:euclidupper} to a setting where one space is a compact subset of a Riemannian manifold and the other has polynomial covering numbers. The rate exponent depends on the manifold dimension, even if the covering exponent of the other space is larger. For a compact metric space $S$, let $N_S(\eps)$ be the internal covering number of $S$ with covering radius $\eps>0$.

\begin{theorem}[Riemannian upper bound]\label{thm:riem}
Let $p >0,\, q\geq 1$, and let $M$ be a complete, connected Riemannian manifold of dimension $d\ge 1$ without boundary, with
geodesic distance $d_M$. Suppose that $\X$ is a compact subset of $M$ with the
restriction of $d_M$. Let $(\Y,\dY)$ be a compact metric space satisfying
$N_\Y(\eps)\le A\eps^{-d'}$ for some $A,d'>0$ and all $0<\eps\le1$.
Then, for all $\mu\in\PP(\X)$, $\nu\in\PP(\Y)$ and $n,m\ge2$,
\begin{equation}\label{eq:upper}
  \E\bigl|D_{p,q}^p(\hat\mu_n,\hat\nu_m)-D_{p,q}^p(\mu,\nu)\bigr|\le C\bigl[\psi_{d,p\wedge2}(n)+\psi_{d,p\wedge2}(m)\bigr],
\end{equation}
where $C$ depends only on $p$, $q$, $A$, $d'$, $M$, $\X$ and $\Y$.
\end{theorem}
If both $\X$ and $\Y$ are compact subsets of complete, connected Riemannian manifolds, Theorem~\ref{thm:riem} can be applied with the smaller of the two manifold dimensions. Indeed, each compact subset has polynomial covering numbers, and the two spaces and their samples can be interchanged. Thus, taking both manifolds to be Euclidean spaces, the upper bound in Theorem~\ref{thm:euclid}(a) follows immediately.

The dependence of $C$ on $M$ comes from the choice of charts and
regularity estimates in these charts.
In the Euclidean case~\eqref{eq:euclidupper}, this dependence is absorbed into the
dependence on the ambient dimensions.

\subsection{Proof overview}\label{sec:overview}
We give an overview of the proofs of our main results. Steps~1--3 prove the Riemannian upper bound~\eqref{eq:upper}
for $p>0$ and $q\ge1$, from which the Euclidean upper
bound~\eqref{eq:euclidupper} follows immediately, and Step~4 proves
the Euclidean lower bound~\eqref{eq:lowerintro} for all $p,q>0$.
By homogeneity of $D_{p,q}$ (see the beginning of Section~\ref{sec:reduction}), we assume in Steps~1--3
that $\diam\X,\diam\Y\le1$. As in Theorem~\ref{thm:riem}, $\X$ is the space whose ambient manifold dimension
$d$ enters the rate.

Two difficulties arise in the proof of the upper bound~\eqref{eq:upper}. First, $D_{p,q}^p(\mu,\nu)$ is the minimum
of a nonconvex quadratic functional of the coupling. Unlike for optimal transport, no general duality theory is known
for $(p,q)$-GW on compact metric spaces. Second, the rate exponent depends only on the dimension $d$ of $\X$, although $\Y$ may have a much
larger covering exponent $d'$. Step~1
addresses the first difficulty and reduces the problem to two empirical processes using local duality, one for each
sample. Step~2 bounds the process for the sample from $\X$; it is there that the index $p\wedge2$ arises. Step~3 controls the contribution of the sample from $\Y$, which is where the second
difficulty lies.
Step~1 holds for all $p,q>0$ (see Remark~\ref{rem:general}), and we use
the assumption $q\ge1$ only in Steps~2 and~3.

\subsubsection*{Step 1: reduction to empirical processes via local duality}

The goal of this step is to show (Proposition~\ref{prop:stab}) that the plug-in error is bounded by
two empirical-process suprema, one for each sample; see~\eqref{eq:reduction} below. We first derive by local duality a one-step bound~\eqref{eq:onestep} with a quadratic remainder,
and then iterate it along the line segment joining $(\mu,\nu)$ and $(\hat\mu_n,\hat\nu_m)$ to make
the remainders vanish.

For $z=(x,y),z'=(x',y')\in\X\times\Y$, set
$K_{p,q}(z,z')=\bigl|\dX(x,x')^q-\dY(y,y')^q\bigr|^p$.
For $\gamma\in\PP(\X\times\Y)$, define
\begin{equation}\label{eq:K}
  J_{p,q}(\gamma)
  =\iint K_{p,q}(z,z')\,d\gamma(z)\,d\gamma(z'),
\end{equation}
so that $D_{p,q}^p(\mu,\nu)
=\inf_{\gamma\in\Pi(\mu,\nu)}J_{p,q}(\gamma)$.

\medskip\noindent\emph{Local duality.}
We expand the quadratic objective $J_{p,q}$ around a minimizer and apply optimal transport
duality to its linear term, which gives the one-step bound~\eqref{eq:onestep}.
To motivate this, recall that for optimal transport with a cost $c$, the reduction to
empirical processes follows directly from Kantorovich duality. Write $T_c(\mu',\nu')$ for the optimal transport value. If $(f,g)$ is an optimal Kantorovich potential pair for $(\hat\mu_n,\hat\nu_m)$, then strong duality at $(\hat\mu_n,\hat\nu_m)$ and weak
duality at $(\mu,\nu)$ give
\begin{equation}\label{eq:OTdual}
  T_c(\hat\mu_n,\hat\nu_m)-T_c(\mu,\nu)\le\int f\,d(\hat\mu_n-\mu)+\int g\,d(\hat\nu_m-\nu).
\end{equation}
Exchanging the roles of the two marginal pairs gives a bound in the reverse direction, using an
optimal pair for $(\mu,\nu)$ that may differ from $(f,g)$ (see, for instance, the proof of~\citep[Theorem~2.2]{hundrieser2024empirical}).

In contrast, $D^p_{p,q}(\mu,\nu)$ has no such duality in general. In the Euclidean case, previous works treating even integers $p$
\citep{zhang2024gromov,leung2026sample,paliy2026empirical} address this issue by
representing $D_{p,q}^p$ variationally as an optimization over a class of
optimal transport problems with parameter-dependent costs, to which Kantorovich
duality applies. This approach does not directly extend to other $p$, since it relies on
expanding the $p$-th power in $K_{p,q}$.

The objective $J_{p,q}$, however, is quadratic in the coupling, and we expand it to recover a
transport problem. Writing
$$c_\gamma(z):=2\int K_{p,q}(z,z')\,d\gamma(z'),$$
we have, for $\gamma,\pi\in\PP(\X\times\Y)$ and $\lambda\in[0,1]$,
$$
  J_{p,q}\bigl(\gamma+\lambda(\pi-\gamma)\bigr)-J_{p,q}(\gamma)
  =\lambda\int c_\gamma\,d(\pi-\gamma)+\lambda^2J_{p,q}(\pi-\gamma),
$$
where $J_{p,q}$ is extended to signed measures by~\eqref{eq:K}. Set
$(\mu_\lambda,\nu_\lambda):=((1-\lambda)\mu+\lambda\hat\mu_n,(1-\lambda)\nu+\lambda\hat\nu_m)$.
Choose $\gamma$ to be optimal for the problem defining $D^{p}_{p,q}(\mu,\nu)$, and choose $\pi$ to be an optimal plan for $T_{c_\gamma}(\hat\mu_n,\hat\nu_m)$. Since $\gamma\in\Pi(\mu,\nu)$, the linear term $\int c_\gamma\,d(\pi-\gamma)$ is at most
$T_{c_\gamma}(\hat\mu_n,\hat\nu_m)-T_{c_\gamma}(\mu,\nu)$, to which~\eqref{eq:OTdual} applies.\footnote{In fact, equality holds, since $\gamma$ is also an optimal plan for $T_{c_\gamma}(\mu,\nu)$, as noted in Dumont et al.~\citep[Section~1.1.2]{dumont2025existence}. Thus $\int c_\gamma\,d(\pi-\gamma)
=T_{c_\gamma}(\hat\mu_n,\hat\nu_m)-T_{c_\gamma}(\mu,\nu)$.} Since $\gamma+\lambda(\pi-\gamma)\in\Pi(\mu_\lambda,\nu_\lambda)$ and $J_{p,q}(\gamma)=D_{p,q}^p(\mu,\nu)$, we obtain, for $\lambda\in[0,1]$ (see~\eqref{eq:local} with $t=0$),
\begin{equation}\label{eq:onestep}
  D_{p,q}^p(\mu_\lambda,\nu_\lambda)-D_{p,q}^p(\mu,\nu)
  \le\lambda\Bigl[\int f\,d(\hat\mu_n-\mu)+\int g\,d(\hat\nu_m-\nu)\Bigr]
  +\lambda^2J_{p,q}(\pi-\gamma),
\end{equation}
where $(f,g)$ is an optimal Kantorovich potential pair for $T_{c_\gamma}(\hat\mu_n,\hat\nu_m)$. We call this use of duality \emph{local duality}.

\medskip\noindent\emph{Iterating along the line segment.}
Taking $\lambda=1$ in~\eqref{eq:onestep} would compare $(\mu,\nu)$ with $(\hat\mu_n,\hat\nu_m)$ directly, but the general bound
$|J_{p,q}(\pi-\gamma)|\le4$ for the quadratic term does
not decay with $n$ and $m$. We therefore move from $(\mu,\nu)$ to
$(\hat\mu_n,\hat\nu_m)$ in many small steps along the line segment $t\mapsto(\mu_t,\nu_t)$, $t\in[0,1]$,
applying~\eqref{eq:onestep} at each step with the current marginals $(\mu_t,\nu_t)$ in place of $(\mu,\nu)$ and a new minimizer
$\gamma_t\in\Pi(\mu_t,\nu_t)$. Over a fine subdivision of $[0,b]$ with $b<1$, the linear
terms add up while the quadratic terms become negligible as the mesh tends to $0$. Letting $b\uparrow1$ and interchanging the marginal pairs then gives
the bound of Proposition~\ref{prop:stab},
\begin{equation}\label{eq:reduction}
  \bigl|D_{p,q}^p(\hat\mu_n,\hat\nu_m)-D_{p,q}^p(\mu,\nu)\bigr|
  \le\sup_{f\in\mathcal F_\X}\Bigl|\int f\,d(\hat\mu_n-\mu)\Bigr|
  +\sup_{g\in\mathcal F_\Y}\Bigl|\int g\,d(\hat\nu_m-\nu)\Bigr|,
\end{equation}
which is precisely the reduction needed. Here, $\mathcal F_\X$ and $\mathcal F_\Y$ are the classes of first and second components of $c_\gamma$-conjugate potential pairs $(f,g)$ normalized by $\inf_\X f=0$, with $\gamma$ ranging over $\PP(\X\times\Y)$. The suprema are needed because the minimizers $\gamma_t$, and
with them the potentials, depend on the samples. Both substeps do not use the specific form of $K_{p,q}$, so that
this approach applies to all $p,q>0$ (Remark~\ref{rem:general}).

The price of this reduction is that, since $c_\gamma$ varies with
$\gamma$, the suprema must be taken over potentials for all couplings. It remains to show that suitable bounds on the metric entropies of the potential classes $\mathcal F_\X$ and $\mathcal F_\Y$ can nevertheless be obtained for the two empirical-process suprema
in~\eqref{eq:reduction} to decay in expectation at the rate $\psi_{d,p\wedge2}$,
through Dudley's entropy bound (Lemma~\ref{lem:dudley}). Step~2 obtains such a bound for
$\mathcal F_\X$, and Step~3 for $\mathcal F_\Y$.

\subsubsection*{Step 2: controlling potentials on $\X$}
The goal of this step is to show that (Proposition~\ref{prop:Xside})
$$
  \E\sup_{f\in\mathcal F_\X}
  \Bigl|\int f\,d(\hat\mu_n-\mu)\Bigr|\le C\psi_{d,p\wedge2}(n).
$$
By Dudley's entropy bound, it would suffice to show that the potentials in $\mathcal F_\X$ form a
class with sup-norm metric entropy of order $\eps^{-d/(p\wedge2)}$. We obtain
such a bound for the potentials themselves when $p\le1$. When $p>1$, we
obtain it for corrected potentials, which differ from those in
$\mathcal F_\X$ by a function of $x$, and control the correction
separately. In both cases, we establish uniform regularity estimates for
compositions of $t_a:=d_M(\cdot,a)^q$, $a\in\X$, with suitable scalar
functions and transfer these estimates to the relevant potentials through
the corresponding $c$-transform.

\medskip\noindent\emph{Reduction to regularity.}
We first reduce the entropy bound to uniform regularity of the functions from which the costs
$c_\gamma$ are built.
Every $f\in\mathcal F_\X$ can be written as the $c_\gamma$-transform
$f(x)=\min_{y\in\Y}\{c_\gamma(x,y)-g(y)\}$ for some
$\gamma\in\PP(\X\times\Y)$ and some $g\in\mathcal F_\Y$.
Writing $\varphi_\ell(t):=|t-\ell|^p$ for $\ell\ge0$, we have
\begin{equation}\label{eq:cost-mixture}
  c_\gamma(x,y)
  =2\int\varphi_{\dY(y,b)^q}\bigl(t_a(x)\bigr)\,d\gamma(a,b).
\end{equation}
Thus, for each $y\in\Y$, the function $c_\gamma(\cdot,y)$ is, up to the factor $2$, a mixture
of the functions $\varphi_\ell\circ t_a$ over $\gamma$ with $a\in\X$ and
$\ell=\dY(y,b)^q\in[0,1]$. Hölder or semiconcavity bounds on these functions, uniform in $a\in\X$ and $\ell\in[0,1]$, would pass through
mixtures and pointwise minima, and hence to the potentials through the
$c_\gamma$-transform. The resulting bounds are uniform in $\gamma$, which is essential since the potentials vary over all couplings. Here semiconcavity means an upper bound $C\|v\|^2$
on the second differences
$\Delta_v^2w(z):=w(z+v)+w(z-v)-2w(z)$ in a chart of $M$, with semiconcavity
constant $C$.\footnote{Since chart images need not lie in $\X$, we extend $c_\gamma$ to $M\times\Y$
and $f$ to $M$ by~\eqref{eq:extended-cost} and~\eqref{eq:extended-potential}, respectively. We restrict the potentials
to $\X$ for the entropy bounds in Lemmas~\ref{lem:entHolder} and~\ref{lem:entE}, and suppress
the extension notation in this overview.}

For $p\le1$, we use $p$-Hölder continuity of the compositions $\varphi_\ell\circ t_a$. For $p>1$, we use
semiconcavity, smoothing $\varphi_\ell$ when $1<p<2$, and
subtracting from the cost a function of $x$ alone for all $p>1$, following the
correction argument in~\citep[Section~5.1, Step~3.1]{leung2026sample}. The
corresponding correction to the potentials contributes only a term of order
$n^{-1/2}$ to the bound.

\medskip\noindent\emph{Charts and distance functions.}
Since $\X$ is compact, it is covered by finitely many charts of $M$,
in which we work throughout. We first record the regularity
of the functions $t_a=d_M(\cdot,a)^q$, $a\in\X$, in these charts. First, these functions are uniformly Lipschitz (Lemma~\ref{lem:riemR}),
and this is all that the case $p\le1$ uses.\footnote{Here the proof uses $q\ge1$, which ensures that $u\mapsto u^q$
is Lipschitz on bounded intervals in $[0,\infty)$ (see
Appendix~\ref{app:riem}).} For $p>1$, we also need one-sided
second-difference bounds for these functions in the charts, namely a lower
bound near $a$ and an upper bound away from $a$,
where they are uniformly semiconcave, with an overlapping
region where both hold. In particular, writing $t(z):=(t_a\circ\psi)(z)=d_M(\psi(z),a)^q$ in a chart
$\psi$, Lemma~\ref{lem:riemR} gives constants $\rho,K_t>0$ such that,
for small $v$,
\begin{equation}\label{eq:onesided}
  \Delta_v^2t(z)\ge-K_t\|v\|^2\ \text{ if } d_M(\psi(z),a)<2\rho,\qquad
  \Delta_v^2t(z)\le K_t\|v\|^2\ \text{ if } d_M(\psi(z),a)>\rho.
\end{equation}

We then derive uniform regularity estimates for the compositions
$\varphi_\ell\circ t_a$ by combining these bounds on $t_a$ with those on
$\varphi_\ell$, using smoothing and correction as needed for $p>1$.
From these estimates we obtain entropy bounds for the potentials and
deduce the target rate $\psi_{d,p\wedge2}(n)$.

\medskip\noindent\emph{The case $p\le1$.}
For $p\le1$, the compositions $\varphi_\ell\circ t_a$ are uniformly $p$-Hölder,
since the functions $t_a$ are uniformly Lipschitz and
$t\mapsto|t-\ell|^p$ is $p$-Hölder uniformly in $\ell$.
The potentials $f\in \mathcal{F}_\mathcal{X}$ inherit this bound, so that they are uniformly bounded and $p$-Hölder in each of the
charts. The standard Hölder entropy estimate \citep[Theorem~2.7.1]{vaart1996weak} therefore bounds
their metric entropy in the sup norm by a multiple of $\eps^{-d/p}$ (Lemma~\ref{lem:entHolder}),
and Dudley's entropy bound gives the desired rate $\psi_{d,p}(n)$. This proves
Proposition~\ref{prop:Xside} for $p\le1$.

\medskip\noindent\emph{The case $p>1$.}
For $p>1$, a $p$-Hölder bound would force a function on a convex
chart domain to be constant, while Lipschitz bounds alone give an
entropy bound of order $\eps^{-d}$, which is too weak for the target
bound $\eps^{-d/(p\wedge2)}$.
We thus turn to semiconcavity, but there are two obstacles to obtaining it uniformly for the compositions $\varphi_\ell\circ t_a$ and hence for the potentials in $\mathcal{F}_\mathcal{X}$.
First, for $1<p<2$, $\varphi_\ell$ is not semiconcave on any neighborhood of $\ell$, where its second derivative is
unbounded; in particular, for $0<\ell<1$ and small $v$, $\Delta_v^2\varphi_\ell(\ell)=2|v|^p$, which is not $O(|v|^2)$. Second, for
every $p>1$, the sign of $\varphi_\ell'$ need not match the one-sided second-difference bounds
for $t_a$ in~\eqref{eq:onesided} (see Remark~\ref{rem:correction}). We remove the first obstacle by smoothing $\varphi_\ell$ and the second
by subtracting a correction independent of $\ell$ before composing
with $t_a$.

\emph{Smoothing.}
For $1<p<2$, we obtain a semiconcave approximation to $\varphi_\ell$ by smoothing. We replace
$\varphi_\ell$ by $\varphi_{\ell,h}(t):=((t-\ell)^2+h^2)^{p/2}$, which satisfies
$0\le\varphi_{\ell,h}-\varphi_\ell\le h^p$ and $\varphi_{\ell,h}''\lesssim h^{p-2}$. Using
$\varphi_{\ell,h}$ in place of $\varphi_\ell$ in the
formula~\eqref{eq:cost-mixture} for $c_\gamma$ turns each potential into a
smoothed potential (see~\eqref{eq:smoothed}), which differs from the
original potential by at most $2h^p$ in the sup norm
(Lemma~\ref{lem:smooth}). For $p\ge2$, $\varphi_\ell$ is already semiconcave on the range of $t_a$ in these charts, so no smoothing is
needed, and we take $h=0$, so that $\varphi_{\ell,0}=\varphi_\ell$.

\emph{Correction.}
We next construct the correction that makes the compositions semiconcave in charts, uniformly in $a$ and $\ell$. With the function $t$ as in~\eqref{eq:onesided}, Taylor's inequality gives
$$
\Delta_v^2(\varphi_{\ell,h}\circ t)(z)
\le \varphi_{\ell,h}'(t(z))\,\Delta_v^2t(z)+C\|v\|^2,
$$
where $C$ depends on an upper bound for $\varphi_{\ell,h}''$ and the
Lipschitz constant of $t$. Thus, the sign of $\varphi_{\ell,h}'$ determines
whether we need an upper or a lower bound on the second differences
of $t$ (Lemma~\ref{lem:comp}) to attain uniform semiconcavity. Combined with~\eqref{eq:onesided}, this means that we want $\varphi_{\ell,h}\circ t$ to be nonincreasing for small $t$, that is, for $x$ near $a$, and
nondecreasing for larger $t$. The functions
$\varphi_{\ell,h}$ do not satisfy these conditions uniformly in $\ell$, and without a correction the
potentials need not be semiconcave.

To enforce them, we fix a constant $L_p>0$, depending only on $p$ and $q$, such that
$|\varphi'_{\ell,h}|\le L_p$ on the range of $t_a$ in these charts for all $\ell$ and $h$. We then subtract
from $\varphi_{\ell,h}$ a single $C^2$ correction function $\Omega$, independent of $\ell$ and $h$, whose
derivative equals $L_p$ for small $t$ and $-L_p$ for larger $t$; the transition between the two values
takes place in the overlap region.
Thus $\varphi_{\ell,h}-\Omega$ has the required monotonicity for small and for larger $t$. For $\gamma\in\PP(\X\times\Y)$, define the
corresponding correction to the cost and the potential by
\[
  \Omega_\gamma(x):=2\int\Omega\bigl(t_a(x)\bigr)\,d\gamma(a,b),
\]
a function of $x$ that depends on $\gamma$ only through its first marginal. For each $h$, the corrected functions $(\varphi_{\ell,h}-\Omega)\circ t_a$
are semiconcave in charts uniformly over $a$ and $\ell$, with a semiconcavity constant that may
depend on $h$ (Lemma~\ref{lem:key}).

Finally, subtracting $\Omega_\gamma$ from $f\in\mathcal F_\X$ and from its smoothed
version gives the corrected potential $f^\circ:=f-\Omega_\gamma$ and the
smoothed corrected potential (see~\eqref{eq:corrected}). The latter is the $c_{\gamma,h}^\circ$-transform of
some $g\in\mathcal F_\Y$, where, as in~\eqref{eq:cost-mixture} for $c_\gamma$, the smoothed corrected
cost $c_{\gamma,h}^\circ(\cdot,y)$ is a mixture of
the functions $(\varphi_{\ell,h}-\Omega)\circ t_a$ over $\gamma$ up to the factor $2$.
Hence, for each fixed $h$, the smoothed corrected potentials are
uniformly semiconcave in charts (Lemma~\ref{lem:key}).

\emph{Entropy and rate.}
It remains to derive entropy bounds from the semiconcavity of the smoothed corrected potentials and apply Dudley's entropy bound.
Starting from the entropy bound for the smoothed corrected potentials, we
pass to the corrected potentials, each of which is within $2h^p$ of its
smoothed version in the sup norm, and finally to the potentials in
$\mathcal F_\X$ by bounding the contribution of $\Omega_\gamma$ to the
empirical process.

The smoothed corrected potentials are uniformly bounded and Lipschitz, and
semiconcave in charts with constant $2\Lambda_h$, where $\Lambda_h\lesssim h^{p-2}$ for $1<p<2$, or $\Lambda_h\lesssim1$
for $p\ge2$, uniformly in $\gamma$ (Lemma~\ref{lem:key}).
A function that is semiconcave with constant $\lambda$ is, in a chart, a concave function
plus $\lambda\|z\|^2/2$. Bronshtein's entropy bound (Theorem~\ref{thm:bron}) for uniformly bounded, uniformly Lipschitz convex
functions (\cite{bronshtein1976entropy}; \cite[Corollary~2.7.10]{vaart1996weak}; see also \cite{guntuboyina2013covering}) therefore applies
to the smoothed corrected potentials in finitely many charts. Their sup-norm metric entropy at radius
$\eps$ is bounded by a multiple of $((1+\Lambda_h)/\eps)^{d/2}$. For $1<p<2$, we choose
$h\asymp\eps^{1/p}$, so that the smoothing error $2h^p$ is a small fraction of $\eps$, and the
resulting entropy bound for the corrected potentials $f^\circ$ is
\[
  \left(\frac{\eps^{(p-2)/p}}{\eps}\right)^{d/2}=\eps^{-d/p}.
\]
For $p\ge2$, we take $h=0$, so there is no smoothing error, and the bound is a multiple of
$\eps^{-d/2}$ (Lemma~\ref{lem:entE}). In both cases, Dudley's entropy bound shows that the expected supremum
of the empirical process indexed by the corrected potentials $f^\circ$
is at most $C\psi_{d,p\wedge2}(n)$.

The correction itself is inexpensive to control. Since $f=f^\circ+\Omega_\gamma$ and $\Omega_\gamma$
depends linearly on the first marginal of $\gamma$, writing $\mathcal E$ for the class of corrected
potentials (see Section~\ref{sec:entropy}), we have
$$
  \E\sup_{f\in\mathcal F_\X}\Bigl|\int f\,d(\hat\mu_n-\mu)\Bigr|
  \le\underbrace{\E\sup_{f^\circ \in\mathcal E}\Bigl|\int f^\circ \,d(\hat\mu_n-\mu)\Bigr|}_{\lesssim\,\psi_{d,p\wedge2}(n)}
  +\underbrace{2\,\E\sup_{a\in\X}\Bigl|\int(\Omega\circ t_a)\,d(\hat\mu_n-\mu)\Bigr|}_{\lesssim\,n^{-1/2}}.
$$
Dudley's bound for the corrected potentials gives the first term. For the second, the map $a\mapsto\Omega\circ t_a$ is
Lipschitz in the sup norm and $\X$ has polynomial covering numbers, so this family has logarithmic
metric entropy (Lemma~\ref{lem:entE}). Hence the correction does not affect the rate, which proves
Proposition~\ref{prop:Xside} for $p>1$.

\medskip
In each regime, the entropy exponent is $d/(p\wedge2)$. For $p\le1$, it comes from Hölder
continuity; for $1<p<2$, from the trade-off between the smoothing error $h^p$ and the
semiconcavity constant $\Lambda_h\lesssim h^{p-2}$; and for $p\ge2$, from semiconcavity of the
corrected potentials. This explains the index $p\wedge2$ in the rate, while $q$ enters only
through the leading constants.

This bounds the first supremum in~\eqref{eq:reduction}. Step~3 bounds the second, reusing the
corrected potentials and their entropy bound from this step.

\subsubsection*{Step 3: controlling potentials on $\Y$}
The goal of this step is to show that (Proposition~\ref{prop:Yside})
$$
  \E\sup_{g\in\mathcal F_\Y}
  \Bigl|\int g\,d(\hat\nu_m-\nu)\Bigr|\le C\psi_{d,p\wedge2}(m).
$$
The argument of Step~2 does not apply directly to $\Y$, since $\Y$ need not
have a Riemannian structure. Even if it does, its dimension may be larger
than $d$. We therefore turn to the $c$-transform used to express
$g\in\mathcal F_\Y$ (see~\eqref{eq:g-ctransform} below). Following the
approach of~\citep{zhang2024gromov,leung2026sample,paliy2026empirical},
which is based on the lower complexity adaptation
principle~\citep{hundrieser2024empirical}, we use this representation to
transfer the entropy bound for the corrected potentials from Step~2 to
$\mathcal F_\Y$. By Dudley's entropy bound (Lemma~\ref{lem:dudley}), it
suffices to estimate the entropy of $\mathcal F_\Y$ in $L^2(Q)$, uniformly
over finitely supported probability measures $Q$ on $\Y$. Since the
corrected costs in the $c$-transform vary with $\gamma$, we must also
control their metric entropy. To this end, we discretize them by replacing the minimum over
$\X$ in the $c$-transform with a minimum over a finite net.
The resulting discretization error is controlled using the uniform Hölder
continuity of the corrected costs and potentials in $x$. It then suffices to cover the
discretized costs in the $L^2(Q)$ norm of their maximum over the net,
which we further control using an $L^r$ norm for a suitable choice of $r\geq 2$.

\medskip\noindent\emph{Reduction to the corrected costs.}
We first construct a cover of $\mathcal F_\Y$ from covers of the corrected potentials and discretized costs.
Let $c_\gamma^\circ:=c_\gamma-\Omega_\gamma$ be the corrected cost, where
$\Omega_\gamma=0$ if $p\le1$. Since $\Omega_\gamma$ depends only on $x$, each $g\in\mathcal F_\Y$ is
the $c_\gamma^\circ$-transform of the corrected potential $f^\circ$ from
Step~2 for a suitable $\gamma\in\PP(\X\times\Y)$, that is,
\begin{equation}\label{eq:g-ctransform}
  g(y)=\min_{x\in\X}\{c_\gamma^\circ(x,y)-f^\circ(x)\}.
\end{equation}
Fix a finitely supported probability measure $Q$ on $\Y$ and a target
covering radius $\eps$. To use the $L^r$ convex-hull entropy bound for the corrected costs, we replace the minimum over
$\X$ by one over an $\eta$-net $\{x_1,\ldots,x_k\}$, denoting the resulting function by $g_\eta$.
Using $q\ge1$, we see that the corrected costs and potentials are
uniformly $(p\wedge1)$-Hölder in $x$, which,
by~\eqref{eq:g-ctransform}, implies that the discretization error in
the sup-norm approximation of $g$ is $O(\eta^{p\wedge1})$. We choose
$\eta\asymp\eps^{1/(p\wedge1)}$ with a sufficiently small constant so that
the error is at most $\eps/3$, using a net of size $k\lesssim\eps^{-d/(p\wedge1)}$.

Let $P=Q\otimes \mathrm{Unif}_k$, where $\mathrm{Unif}_k$ is uniform on $\{1,\ldots,k\}$.
Here the discretized corrected costs are the functions $H_\gamma(y,i):=c_\gamma^\circ(x_i,y)$ on
$\Y\times\{1,\ldots,k\}$. The elementary inequality $|\min_i u_i-\min_i v_i|\le\max_i|u_i-v_i|$
and the comparison of norms give, for $r\ge2$,
$$
  \|g_\eta-g'_\eta\|_{L^2(Q)}
  \le \|f^\circ-f'^\circ\|_{\infty,\X}
  +k^{1/r}\|H_\gamma-H_{\gamma'}\|_{L^r(P)},
$$
where $g'_\eta$ is defined analogously from $\gamma'$ and $f'^\circ$.
Thus, covers of the corrected potentials in the sup norm and of the
discretized corrected costs in $L^r(P)$ yield a cover of $\mathcal F_\Y$ in $L^2(Q)$, after allowing
for the factor $k^{1/r}$ and the discretization error (Lemma~\ref{lem:discrete}).

\medskip\noindent\emph{Entropy of the corrected costs.}
We next bound the entropy of the discretized corrected costs by a convex-hull argument adapted
from~\citep[Lemma~5.5]{leung2026sample}.
Each $H_\gamma$ is a mixture over $\gamma$ of the maps
$H_\theta:(y,i)\mapsto c_\theta^\circ(x_i,y)$, $\theta=(a,b)\in\X\times\Y$, where $c_\theta^\circ:=c_{\delta_\theta}^\circ$ is
the corrected cost for the point mass $\delta_\theta$. As mixtures of these uniformly bounded maps,
the discretized corrected costs belong to the closed convex hull of
$\{H_\theta:\theta\in\X\times\Y\}$ in $L^r(P)$. Since $\X$ and $\Y$ have polynomial covering
numbers and $H_\theta$ depends Hölder continuously on $\theta$, this generating class also has
polynomial covering numbers, with constants independent of $Q$, $k$, and $\eps$.
The $L^r$ convex-hull entropy theorem of
Carl et al.~\citep{carl1999metric} (see also \citep{song2002upper}) therefore applies to the discretized corrected costs, bounding
their entropy at radius $\eps$ by a multiple of
$\eps^{-\kappa}$, where $\kappa<2$.
To account for the factor $k^{1/r}$, we cover the discretized corrected costs with radius
$\eps/(3k^{1/r})$ in $L^r(P)$. Since $k\lesssim\eps^{-d/(p\wedge1)}$, their metric
entropy at this radius is $O(\eps^{-\kappa_r})$, where
$\kappa_r=\kappa(1+d/((p\wedge1)r))$. Since $\kappa<2$ is independent of $r$, we can fix $r$ large
enough that $\kappa_r<2$ (Lemma~\ref{lem:costs}).

\medskip\noindent\emph{Entropy and rate.}
Finally, we combine the two covers and apply Dudley's entropy bound.
With an $\eps/3$-cover of the discretized corrected costs and a
sup-norm $\eps/3$-cover of the corrected potentials
from~\eqref{eq:Eentropy}, the two covering errors and the
discretization error sum to at most $\eps$. This gives a metric entropy bound of order
$$
  \eps^{-d/(p\wedge2)}+\eps^{-\kappa_r}
$$
for $\mathcal F_\Y$ in $L^2(Q)$, uniformly in $Q$ (Lemma~\ref{lem:entY}). The second term contributes only $O(m^{-1/2})$
through Dudley's entropy bound, while the first gives
$\psi_{d,p\wedge2}(m)$. This proves Proposition~\ref{prop:Yside}.

Together with Step~2 and the reduction~\eqref{eq:reduction}, this proves~\eqref{eq:upper}.

\subsubsection*{Step 4: lower bound}
We now outline the proof of the lower bound~\eqref{eq:lowerintro}, which is independent of the upper-bound argument.
By isometric embedding, it suffices to construct two measures $\mu$ and $\nu$ on $B_1^d(0)$ with $d:=d_\star$ that attain the rates of the upper bound~\eqref{eq:euclidupper}, up to the logarithmic factor in the critical dimension $d = 2(p\wedge2)$. We prove the cases $d<2(p\wedge2)$ and $d \geq  2(p\wedge2)$ separately.

\medskip\noindent\emph{The case $d<2(p\wedge2)$.}
Proposition~\ref{prop:lower-small-dim} gives the lower bound $c(m\wedge n)^{-1/2}$ in~\eqref{eq:lowerintro}, using measures $\mu = (1-\theta) \delta_0 + \theta \delta_{e_1}$ and $\nu = \delta_0$, where $e_1$ is the first standard coordinate vector and $\theta \in (0,1) \setminus \{1/2\}$. We show that $D_{p,q}^p(\mu,\nu) = 2\theta (1-\theta)$ using these measures and bound $\E\bigl|D_{p,q}^p(\hat\mu_n,\hat\nu_m)-D_{p,q}^p(\mu,\nu)\bigr|$ using the mean absolute deviation of $\hat\theta_n$ from $\theta$, where $\hat\mu_n = (1-\hat \theta_n) \delta_0 + \hat\theta_n \delta_{e_1}$.

\begin{figure}[!htbp]
\centering
\definecolor{RhoPoint}{HTML}{29313D}
\definecolor{RhoCube}{HTML}{157F88}
\definecolor{RhoAxis}{HTML}{7F8790}
\resizebox{0.6\linewidth}{!}{%
\begin{tikzpicture}[
  font=\normalfont\normalsize,
  axis/.style={draw=RhoAxis,line width=.55pt,-{Stealth[length=4pt,width=3pt]}},
  point/.style={fill=RhoPoint,draw=white,line width=.45pt},
  point label/.style={align=center,inner sep=2pt,text=RhoPoint},
  q label/.style={align=left,inner sep=2pt,text=RhoCube},
  projection/.style={draw=RhoAxis!80,line width=.45pt,dash pattern=on .8pt off 1pt},
  tick/.style={draw=RhoAxis!50,line width=.55pt},
  tick label/.style={font=\small,inner sep=2pt,text=RhoPoint!50},
]
  \draw[draw=RhoAxis!45,line width=.4pt] (7.35,-2.65) -- (7.35,6.3);
  \node at (2.4,6.0) {$d=2$};
  \begin{scope}[x={(6cm,0cm)},y={(0cm,6cm)}]
    \draw[axis] (0,0) -- (.89,0) node[right,text=RhoPoint] {$x_1$};
    \draw[axis] (0,0) -- (0,.89) node[above,text=RhoPoint] {$x_2$};
    \foreach \t in {.125,.25,.375,.5,.625,.75} {
      \draw[tick] ([yshift=-2.5pt]\t,0) -- ([yshift=2.5pt]\t,0);
      \draw[tick] ([xshift=-2.5pt]0,\t) -- ([xshift=2.5pt]0,\t);
    }
    \foreach \t/\lab in {.25/{1/4},.5/{1/2},.75/{3/4}} {
      \node[tick label,anchor=north] at ([yshift=-4pt]\t,0) {$\lab$};
      \node[tick label,anchor=east] at ([xshift=-4pt]0,\t) {$\lab$};
    }
    \pgfmathsetmacro{\lo}{1/(8*sqrt(2))}
    \pgfmathsetmacro{\hi}{3/(16*sqrt(2))}
    \pgfmathsetmacro{\mid}{(\lo+\hi)/2}
    \filldraw[fill=RhoCube!22,draw=RhoCube,line width=.8pt]
      (\lo,\lo) rectangle (\hi,\hi);
    \fill[RhoCube] (\mid,\mid) circle[radius=1pt];
    \node[q label,anchor=south west] at (\hi,\hi) {$Q$};
    \filldraw[point] (0,0) circle[radius=3.1pt];
    \filldraw[point] (5/8,0) circle[radius=2.4pt];
    \filldraw[point] (0,3/4) circle[radius=2.4pt];
    \node[point label,anchor=north east] at (-.025,-.035)
      {$a_0$};
    \node[point label,anchor=south] at (5/8,.045)
      {$a_1$};
    \node[point label,anchor=west] at (.05,3/4)
      {$a_2$};
  \end{scope}
  \node at (12.0,6.0) {$d=3$};
  \begin{scope}[shift={(11.0,0)},
    x={(6cm,0cm)},y={(-2.6cm,-2.6cm)},z={(0cm,6cm)}]
    \draw[axis] (0,0,0) -- (.89,0,0) node[right,text=RhoPoint] {$x_1$};
    \draw[axis] (0,0,0) -- (0,.89,0) node[below left,text=RhoPoint] {$x_2$};
    \draw[axis] (0,0,0) -- (0,0,.89) node[above,text=RhoPoint] {$x_3$};
    \foreach \t in {.125,.25,.375,.5,.625,.75} {
      \draw[tick] ([yshift=-2.5pt]\t,0,0) -- ([yshift=2.5pt]\t,0,0);
      \draw[tick] ([xshift=-1.8pt,yshift=1.8pt]0,\t,0)
        -- ([xshift=1.8pt,yshift=-1.8pt]0,\t,0);
      \draw[tick] ([xshift=-2.5pt]0,0,\t) -- ([xshift=2.5pt]0,0,\t);
    }
    \foreach \t/\lab in {.25/{1/4},.5/{1/2},.75/{3/4}} {
      \node[tick label,anchor=north] at ([yshift=-4pt]\t,0,0) {$\lab$};
      \node[tick label,anchor=south east]
        at ([xshift=-3pt,yshift=3pt]0,\t,0) {$\lab$};
      \node[tick label,anchor=east] at ([xshift=-4pt]0,0,\t) {$\lab$};
    }
    \pgfmathsetmacro{\lo}{1/(8*sqrt(3))}
    \pgfmathsetmacro{\hi}{3/(16*sqrt(3))}
    \pgfmathsetmacro{\mid}{(\lo+\hi)/2}
    \draw[draw=RhoCube!55,densely dashed,line width=.4pt]
      (\hi,\lo,\lo) -- (\lo,\lo,\lo) -- (\lo,\hi,\lo)
      (\lo,\lo,\lo) -- (\lo,\lo,\hi);
    \filldraw[fill=RhoCube!15,draw=RhoCube,line width=.6pt]
      (\lo,\lo,\hi) -- (\hi,\lo,\hi) -- (\hi,\hi,\hi)
      -- (\lo,\hi,\hi) -- cycle;
    \filldraw[fill=RhoCube!35,draw=RhoCube,line width=.6pt]
      (\hi,\lo,\lo) -- (\hi,\hi,\lo) -- (\hi,\hi,\hi)
      -- (\hi,\lo,\hi) -- cycle;
    \filldraw[fill=RhoCube!24,draw=RhoCube,line width=.6pt]
      (\lo,\hi,\lo) -- (\hi,\hi,\lo) -- (\hi,\hi,\hi)
      -- (\lo,\hi,\hi) -- cycle;
    \fill[RhoCube] (\mid,\mid,\mid) circle[radius=1pt];
    \coordinate (qtop) at (\hi,\lo,\hi);
    \filldraw[point] (0,0,0) circle[radius=3.1pt];
    \filldraw[point] (7/12,0,0) circle[radius=2.4pt];
    \filldraw[point] (0,2/3,0) circle[radius=2.4pt];
    \filldraw[point] (0,0,3/4) circle[radius=2.4pt];
    \node[point label,anchor=south east] at (-.035,0,.025)
      {$a_0$};
    \node[point label,anchor=south] at (7/12,0,.045)
      {$a_1$};
    \node[point label,anchor=west] at (.06,2/3,0)
      {$a_2$};
    \node[point label,anchor=west] at (.05,0,3/4)
      {$a_3$};
  \end{scope}
  \node[q label,anchor=south west] at ([yshift=1pt]qtop)
    {$Q$};
\end{tikzpicture}
}
\caption{The measure $\mu$ for $d=2$ (left) and $d=3$ (right).
The points $a_0=0$ and $a_k=(\frac12+\frac{k}{4d})e_k$,
$1\leq k\leq d$, carry probability $2/3$ and $1/(6d)$,
respectively. The remaining probability $1/6$ is uniformly
distributed on the shaded cube
$Q=[1/(8\sqrt d),3/(16\sqrt d)]^d$.
Axis ticks are spaced at intervals of $1/8$.}
\label{fig:rho}
\end{figure}
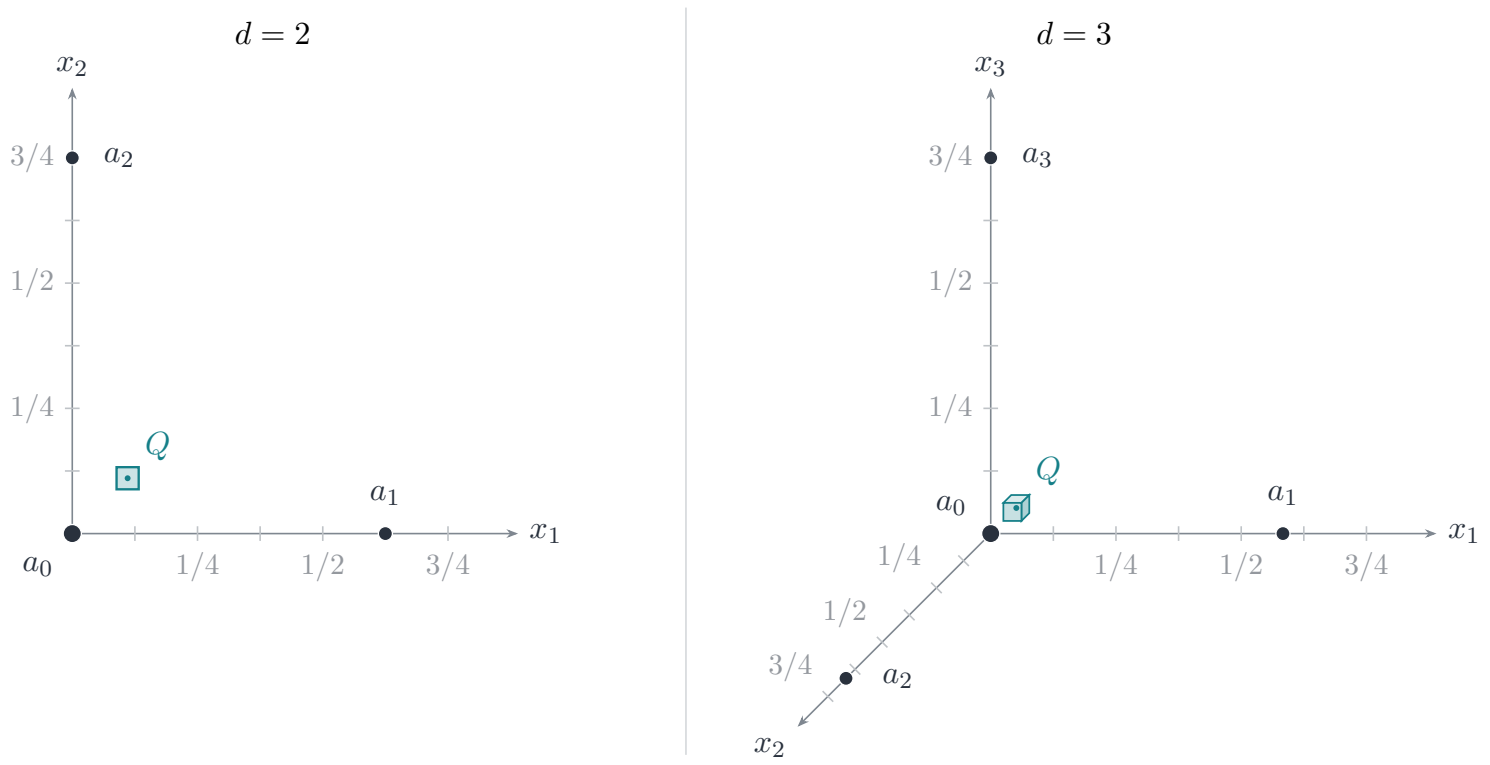

\medskip\noindent\emph{Construction for $d\ge2(p\wedge2)$.}
For $d\ge2(p\wedge2)$, we reduce the problem of lower bounding $\E\bigl|D_{p,q}^p(\hat\mu_n,\hat\nu_m)-D_{p,q}^p(\mu,\nu)\bigr|$ to the estimation of Wasserstein distances on a cube. To construct the measures, let $e_1,\dots,e_d \in \R^d$ be the standard coordinate vectors in $\R^d$. Let $r_0=0$ and $r_k = \frac{1}{2} + \frac{k}{4d}$ for $k = 1,\dots,d$, and let $a_0 = 0$, $a_k := r_ke_k$. Also define the cube $Q := \prod_{i=1}^d \left[\frac{1}{8\sqrt{d}},\frac{3}{16\sqrt{d}} \right]$. Then, define the measures:
\[
    \mu = \frac{2}{3} \delta_{a_0} + \frac{1}{6d} \sum_{k=1}^d \delta_{a_k} + \frac{1}{6} \lambda_Q, \qquad \nu = S_\# \mu, \qquad S(x) = \begin{cases}
        x, \quad &0<p<2 \\
        x/2, \quad &p\geq 2
    \end{cases},
\]
where $\lambda_Q$ is the uniform (i.e. Lebesgue) probability measure on the cube $Q$.
Figure~\ref{fig:rho} illustrates $\mu$ for $d=2$ and $d=3$.

We now explain the construction of $\mu$ and $\nu$ (Section~\ref{app:lower-bound} and Lemma~\ref{lem:fixed-points}). The numerical constants in this construction are chosen for concreteness. To prove their $n^{-2/d}$ lower bound, Zhang et al.~\citep{zhang2024gromov} use measures $\mu = \operatorname{Unif}(B_1^d(0))$ and $\nu = \operatorname{Unif}(B_2^d(0)) = (2\operatorname{id})_\# \mu$. We mirror this by considering uniform measures on a cube $Q$.

\medskip\noindent\emph{Role of the point masses.}
To explain the presence of delta masses, first recall that the GW distance is invariant to isometries, so $D_{p,q}^p(\mu, U_\# \mu) = 0$ for every $U \in E(d)$, where $E(d)$ is the isometry group on $\R^d$. However, $W_s(\mu,U_\# \mu)$ need not be small. A bound of the form $D_{p,q}^p (\mu,\nu)\gtrsim W_s^s(\mu,\nu)$ can therefore only hold modulo isometries, which is why Zhang et al.~\citep{zhang2024gromov} prove a lower bound of the form $D_{2,2}(\mu,\nu)\;\ge\;c_{\mu,\nu}\,\inf_{U\in E(d)}W_2(\mu,U_\#\nu)$, where the infimum over $E(d)$ takes care of the isometry invariance and $c_{\mu,\nu}>0$ depends only on the smallest eigenvalues $\lambda_{\min}(\Sigma_\mu)$ and $\lambda_{\min}(\Sigma_\nu)$ of the covariance matrices of $\mu$ and $\nu$. We address the isometry invariance by placing $d+1$ delta masses $\delta_{a_0},\dots,\delta_{a_d}$ at each of the coordinate directions. We show in the proof of Proposition~\ref{prop:lower-large-dim} that a coupling $\gamma \in \Pi(\mu,\nu)$ that nearly preserves pairwise distances must send each delta mass $\delta_{a_i}$ to its image $\delta_{S(a_i)}$ under $S$ with nontrivial probability. Since an isometry of $\mathbb R^d$ fixing $d+1$ affinely independent points is the identity, we avoid needing to take an infimum over $E(d)$ when lower bounding $D_{p,q}^p(\mu,\nu)$ by a Wasserstein distance.

\medskip\noindent\emph{Comparison with transport costs.}
To proceed with the proof, we consider the cases $0<p < 2$ and $p \geq 2$ separately.

If $0<p<2$, we note that $D_{p,q}^p(\mu,\nu) = 0$ since $\nu = \mu$ by construction, and hence $\E\bigl|D_{p,q}^p(\hat\mu_n,\hat\nu_m)-D_{p,q}^p(\mu,\nu)\bigr|$ equals $\E D_{p,q}^p(\hat\mu_n,\hat\nu_m)$. Next, writing $\hat\nu_m':=(S^{-1})_\#\hat\nu_m$, we restrict the empirical measures to the cube $Q$ by first defining $P_{\mu,n}:=\frac{\hat\mu_n|_Q}{\hat\mu_n(Q)}$ and $P_{\nu,m}:=\frac{\hat\nu_m'|_Q}{\hat\nu_m'(Q)}$, then proving that $D_{p,q}^p(\hat\mu_n,\hat\nu_m) \geq cW_p^p(P_{\mu,n},P_{\nu,m})$ for some $c>0$, on an event $\mathcal E_{n,m}$ of nontrivial probability.

If $p\ge2$, we loosely follow the lower bound proof of Zhang et al.~\citep{zhang2024gromov}. By convexity of $t \mapsto |t|^p$, we lower bound the plug-in error by a difference of two terms (Lemma~\ref{lem:scalar}). The first term depends only on the marginals and has expectation of order $O((m \wedge n)^{-1})$, which is negligible for sufficiently large $m$ and $n$, and the second term may be viewed as a weighted version of $D_{2,q}^2$ and can be lower bounded in expectation by a constant multiple of $\E[\mathbf1_{\mathcal E_{n,m}}W_2^2(P_{\mu,n},P_{\nu,m})]$ (Lemma~\ref{lem:fixed-points}).

\medskip\noindent\emph{Completion of the lower bound.}
Combining these two cases, we conclude the $d \geq 2(p\wedge2)$ case using a simple volumetric argument to lower bound the Wasserstein distances by proving that \[ \E\!\left[\mathbf1_{\mathcal E_{n,m}}W_s^s(P_{\mu,n},P_{\nu,m})\right] \geq c_{s,d}(m \wedge n)^{-s/d}\mathbb P(\mathcal E_{n,m}) \] for $s>0$, where $\mathbb P(\mathcal E_{n,m})\ge1/2$ for sufficiently large $m,n$ (Lemma~\ref{lem:finite-support}). Finally, we combine our lower bounds for the $d < 2(p\wedge2)$ and $d \geq 2(p\wedge2)$ cases to obtain the combined lower bound~\eqref{eq:lowerintro}.
\par

\subsection{Organization}
The rest of the paper is organized as follows. In Sections~\ref{sec:reduction}
and~\ref{sec:thm2a}, we prove the upper bound~\eqref{eq:upper}.
Section~\ref{sec:reduction} establishes the reduction to empirical processes
outlined in Step~1 of the proof overview in Section~\ref{sec:overview}.
Section~\ref{sec:thm2a} completes the proof by bounding the empirical-process suprema, as outlined in Steps~2 and~3.
Section~\ref{app:lower-bound} proves the lower bound~\eqref{eq:lowerintro}.
Appendix~\ref{app:tools} collects auxiliary deterministic lemmas and a version of
Dudley's entropy bound, and Appendix~\ref{app:riem} proves Lemma~\ref{lem:riemR}.

\section{Reduction to empirical processes}\label{sec:reduction}

This section bounds the plug-in error $\bigl|D_{p,q}^p(\hat\mu_n,\hat\nu_m)-D_{p,q}^p(\mu,\nu)\bigr|$ using Kantorovich potentials and local duality, as outlined in Step~1 of the proof overview.

We first normalize the spaces. Since $D_{p,q}^p$ is homogeneous of degree $pq$ in the pair of
metrics, dividing both metrics by $R_1:=\diam\X\vee\diam\Y\vee1$ reduces~\eqref{eq:upper} to the case
$\diam\X,\diam\Y\le1$, with an extra factor $R_1^{pq}$ in the constant. In the setting of
Theorem~\ref{thm:riem}, this multiplies the Riemannian metric on $M$ by $R_1^{-2}$ and preserves the
covering hypothesis on $\Y$, since $N_\Y(R_1\eps)\le N_\Y(\eps)$. After this normalization,
$0\le K_{p,q}\le1$.

For a continuous cost $c$ on $\X\times\Y$, write $T_c(\mu,\nu):=\inf_{\pi\in\Pi(\mu,\nu)}\int c\,d\pi$,
and let $\Phi_{\X,\Y}(c)$ be the set of $c$-conjugate pairs, that is, of $(f,g)\in C(\X)\times C(\Y)$
with
$$
f(x)=\min_{y\in\Y}\bigl(c(x,y)-g(y)\bigr),\qquad
g(y)=\min_{x\in\X}\bigl(c(x,y)-f(x)\bigr).
$$
Such pairs satisfy $f(x)+g(y)\le c(x,y)$ on $\X\times\Y$. For $\mu\in\PP(\X)$ and $\nu\in\PP(\Y)$,
Kantorovich duality \citep[Theorem~1.39]{santambrogio2015optimal} (see also
\citep[Theorem~5.10]{villani2009optimal}) shows that the infimum defining $T_c$ is attained and
$$
  T_c(\mu,\nu)=\max_{(f,g)\in\Phi_{\X,\Y}(c)}\left(\int f\,d\mu+\int g\,d\nu\right).
$$
Since $(f-s,g+s)$ is again $c$-conjugate for every constant $s$, the maximum is attained by a pair
with $\inf_\X f=0$.

For $\gamma\in\PP(\X\times\Y)$, let $c_\gamma(z):=2\int K_{p,q}(z,z')\,d\gamma(z')$, a continuous
function on $\X\times\Y$ with values in $[0,2]$. The following proposition bounds the difference of
$D_{p,q}^p$ between any two marginal pairs by empirical-process suprema. The plug-in error is the
case $(\mu_0,\nu_0)=(\mu,\nu)$, $(\mu_1,\nu_1)=(\hat\mu_n,\hat\nu_m)$.

\begin{proposition}\label{prop:stab}
Let $\mathcal T$ be the set of triples $(\gamma,f,g)$ with $\gamma\in\PP(\X\times\Y)$,
$(f,g)\in\Phi_{\X,\Y}(c_\gamma)$ and $\inf_\X f=0$, and let
$\mathcal F_\X:=\{f:(\gamma,f,g)\in\mathcal T\}$ and $\mathcal F_\Y:=\{g:(\gamma,f,g)\in\mathcal T\}$.
Then, for $(\mu_i,\nu_i)\in\PP(\X)\times\PP(\Y)$, $i=0,1$,
\begin{equation}\label{eq:marginal-stability}
  \bigl|D_{p,q}^p(\mu_1,\nu_1)-D_{p,q}^p(\mu_0,\nu_0)\bigr|
  \le\sup_{f\in\mathcal F_\X}\Bigl|\int f\,d(\mu_1-\mu_0)\Bigr|
  +\sup_{g\in\mathcal F_\Y}\Bigl|\int g\,d(\nu_1-\nu_0)\Bigr|.
\end{equation}
\end{proposition}

\begin{proof}
We first bound the increment of $D_{p,q}^p$ between nearby points on the line segment joining the two
pairs of marginals by local duality, and then sum these one-step bounds along the line segment. The
quadratic remainders vanish as the subdivision is refined.
Write $\mu_t:=(1-t)\mu_0+t\mu_1$, $\nu_t:=(1-t)\nu_0+t\nu_1$ and
$V_t:=D_{p,q}^p(\mu_t,\nu_t)$ for $t\in[0,1]$. Let $B$ denote the sum of the two suprema in~\eqref{eq:marginal-stability}.
For each $t$, a minimizer of $J_{p,q}$ over $\Pi(\mu_t,\nu_t)$ exists, because the set of couplings is
weakly compact and $J_{p,q}$ is weakly continuous, as $K_{p,q}$ is bounded and continuous.

\medskip\noindent\emph{Local duality.}
We first bound the increment of $V$ over a single step.
Let $0\le t<u\le1$, let $\gamma$ minimize $J_{p,q}$ over
$\Pi(\mu_t,\nu_t)$, and set $\lambda:=(u-t)/(1-t)\in(0,1]$. By Kantorovich duality there are an optimal
plan $\pi\in\Pi(\mu_1,\nu_1)$ and an optimal pair $(f,g)$ for the cost
$c_\gamma$, with $(\gamma,f,g)\in\mathcal T$, so that $\int c_\gamma\,d\pi=\int f\,d\mu_1+\int g\,d\nu_1$.
Extending $J_{p,q}$ to finite signed measures on $\X\times\Y$ by~\eqref{eq:K}, for every finite signed measure $\eta$ on $\X\times\Y$, expanding $J_{p,q}(\gamma + \eta)$ using the symmetry of $K_{p,q}$ gives
\begin{equation}\label{eq:expand}
  J_{p,q}(\gamma+\eta)=J_{p,q}(\gamma)+\int c_\gamma\,d\eta+J_{p,q}(\eta).
\end{equation}
With $\eta=\lambda(\pi-\gamma)$,
\[
  J_{p,q}\bigl(\gamma+\lambda(\pi-\gamma)\bigr)-J_{p,q}(\gamma)
  =\lambda\int c_\gamma\,d(\pi-\gamma)+\lambda^2J_{p,q}(\pi-\gamma).
\]
The interpolated plan $\gamma+\lambda(\pi-\gamma)$ couples $(\mu_u,\nu_u)$ and $J_{p,q}(\gamma)=V_t$. Since
$f(x)+g(y)\le c_\gamma(x,y)$, it follows that
\begin{equation}\label{eq:local}
  V_u-V_t\le\lambda\Bigl[\int f\,d(\mu_1-\mu_t)+\int g\,d(\nu_1-\nu_t)\Bigr]
  +\lambda^2J_{p,q}(\pi-\gamma).
\end{equation}
Moreover, $|J_{p,q}(\pi-\gamma)|\le4$, since $0\le K_{p,q}\le1$ and $\pi-\gamma$ has total variation at most $2$. Since $f\in\mathcal F_\X$ and $g\in\mathcal F_\Y$, using $\mu_1-\mu_t=(1-t)(\mu_1-\mu_0)$ and the analogous identity for $\nu_t$ gives
\begin{equation}\label{eq:incr}
  V_u-V_t\le(u-t)B+\frac{4(u-t)^2}{(1-t)^2}.
\end{equation}

\medskip\noindent\emph{Iterating along the line segment.}
We now sum the one-step bounds~\eqref{eq:incr} and pass to the limit.
The factor $(1-t)^{-2}$ in~\eqref{eq:incr} is unbounded as $t\uparrow1$, so we first sum over a
subdivision of $[0,b]$ with $b<1$.
Fix $b\in(0,1)$ and let $t_i:=ib/k$ for $i=0,\ldots,k$. Summing~\eqref{eq:incr} gives
\[
V_b-V_0 =\sum_{i=0}^{k-1}(V_{t_{i+1}}-V_{t_i}) \le bB+\frac{4b^2}{k^2}\sum_{i=0}^{k-1}\frac{1}{(1-t_i)^2}
\le bB+\frac{4b^2}{k(1-b)^2}.
\]
Letting $k\to\infty$ yields
$V_b-V_0\le bB$. To pass to $b=1$, choose $b_j\uparrow1$ along a sequence realizing
$\liminf_{b\uparrow1}V_b$, and choose a minimizer $\gamma_j\in\Pi(\mu_{b_j},\nu_{b_j})$. Compactness
of $\PP(\X\times\Y)$ gives a weakly convergent subsequence whose limit has marginals
$(\mu_1,\nu_1)$. By continuity of $J_{p,q}$, $V_1\le\liminf_{b\uparrow1}V_b$, so
$V_1-V_0\le B$. Interchanging the marginal pairs leaves $B$ unchanged and yields $V_0-V_1\le B$.
\end{proof}

\begin{remark}\label{rem:general}
The proof of Proposition~\ref{prop:stab} uses only that $\X$ and $\Y$
are compact and that $K_{p,q}$ is continuous and symmetric.
It requires neither a Riemannian structure on $\X$ nor a polynomial
covering bound for $\Y$. In particular, the proposition holds for
all $p,q>0$ on arbitrary compact metric spaces. More
generally, it holds with $K_{p,q}$ replaced by any continuous symmetric function $K$ on
$(\X\times\Y)^2$, with $c_\gamma$, $\mathcal T$, $\mathcal F_\X$ and $\mathcal F_\Y$ defined from $K$ and
the bound $|J_{p,q}(\pi-\gamma)|\le4$ replaced by $4\sup|K|$. This covers other variants of GW, for
instance inner-product GW \citep{dandapanthula2025optimal} on compact subsets of Hilbert spaces.
\end{remark}

\section{Proof of the upper bound}\label{sec:thm2a}

We now prove the upper bound~\eqref{eq:upper} by bounding in expectation
the empirical-process suprema in Proposition~\ref{prop:stab}, using metric
entropy estimates for $\mathcal F_\X$ and $\mathcal F_\Y$. Throughout, the distances are
normalized as at the beginning of Section~\ref{sec:reduction}, so that $\diam\X\le1$ and
$\diam\Y\le1$. Constants $C$ may change from line to line and depend only on $p,q,A,d'$ and the
normalized compact subset $\X$ of $M$.

\begin{proposition}[Empirical-process bounds]\label{prop:potentials}
For $n,m\ge2$,
\begin{enumerate}[label=\textup{(\roman*)},ref=\thetheorem(\roman*),nosep]
\item\label{prop:Xside}
\begin{minipage}{\linewidth}
\[
  \E\sup_{f\in\mathcal F_\X}\Bigl|\int f\,d(\muh-\mu)\Bigr|\le C\psi_{d,p\wedge2}(n);
\]
\end{minipage}
\item\label{prop:Yside}
\begin{minipage}{\linewidth}
\[
  \E\sup_{g\in\mathcal F_\Y}\Bigl|\int g\,d(\hat\nu_m-\nu)\Bigr|\le C\psi_{d,p\wedge2}(m).
\]
\end{minipage}
\end{enumerate}
\end{proposition}

Sections~\ref{sec:notation}--\ref{sec:entropy} prove \textup{(i)} by controlling the potentials in
$\mathcal F_\X$, and Section~\ref{sec:yside} proves \textup{(ii)} by controlling those in $\mathcal F_\Y$,
following Steps~2 and~3 of the proof overview. Section~\ref{sec:proof} combines these bounds with
Proposition~\ref{prop:stab} to prove~\eqref{eq:upper}.

\subsection{Notation and reduction to regularity}\label{sec:notation}

We begin the proof of Proposition~\ref{prop:Xside} by fixing notation and
reducing the regularity estimates for the potentials in $\mathcal F_\X$
to estimates for the functions $\varphi_\ell\circ t_a$ defined below.

\medskip\noindent\emph{Sup norms and Hölder continuity.}
For a set $S$ and a bounded function $w$ on $S$, write $\|w\|_{\infty,S}:=\sup_S|w|$, and
$\|w\|_\infty$ if $S$ is clear from context. For a metric space $(S,d_S)$, $C\ge0$ and $s\in(0,1]$, a
function $w$ on $S$ is $Cd_S^s$-Hölder if $|w(x)-w(x')|\le Cd_S(x,x')^s$ for all
$x,x'\in S$, and $s$-Hölder if this holds for some $C$. On subsets of $\R^d$, we use the
Euclidean distance and write $C\|z-z'\|^s$-Hölder.

\medskip\noindent\emph{Covers and covering numbers.}
For a subset $G$ of a seminormed space, $N(\eps,G,\|\cdot\|)$ is the internal covering number of $G$
with covering radius $\eps>0$. Covers are internal unless stated otherwise, and passing from an
external cover to an internal one at most doubles the radius.

\medskip\noindent\emph{Balls and semiconcavity.}
Let $B_r:=\{z\in\R^d:\|z\|<r\}$ be the open Euclidean ball and $\bar B_r$ its closure. For a
real-valued function $w$ on a set $Z\subseteq\R^d$ and $z,z\pm v\in Z$, the second
difference of $w$ is $\Delta_v^2w(z):=w(z+v)+w(z-v)-2w(z)$. A continuous function $w$ on a convex
set $Z\subseteq\R^d$ is \emph{$\lambda$-semiconcave} if $w-\lambda\|\cdot\|^2/2$ is concave,
equivalently, if $\Delta_v^2w(z)\le\lambda\|v\|^2$ whenever $z\pm v\in Z$
(see~\citep[Definition~1.1.1 and Proposition~1.1.3]{cannarsa2004semiconcave}). It suffices to check
this inequality for $\|v\|\le\delta$, for some $\delta>0$, since concavity can be checked locally on a
convex set.

\medskip\noindent\emph{Reduction to regularity.}
We now express the potentials in $\mathcal F_\X$ through the functions $\varphi_\ell\circ t_a$.
For $\ell\ge0$ and $t\ge0$, let $\varphi_\ell(t):=|t-\ell|^p$, and for $a\in\X$ and $x\in M$, let
$t_a(x):=d_M(x,a)^q$. For $\gamma\in\PP(\X\times\Y)$, we extend $c_\gamma$ to $M\times\Y$ by
\begin{equation}\label{eq:extended-cost}
  c_\gamma(x,y):=2\int\varphi_{\dY(y,b)^q}\bigl(t_a(x)\bigr)\,d\gamma(a,b),
\end{equation}
which agrees on $\X\times\Y$ with the definition in Section~\ref{sec:reduction}. For
$(\gamma,f,g)\in\mathcal T$, we extend $f$ to $M$, since the charts of Section~\ref{sec:reg} need
not map into $\X$, by
\begin{equation}\label{eq:extended-potential}
  \bar f(x):=\min_{y\in\Y}\{c_\gamma(x,y)-g(y)\},\qquad x\in M.
\end{equation}

By the conjugacy relations, $\bar f|_\X=f$ and $g(y)=\min_{x\in\X}\{c_\gamma(x,y)-f(x)\}$, and by
Lemma~\ref{lem:pot}, $|f|\le2$ and $|g|\le4$. In view of~\eqref{eq:extended-cost}
and~\eqref{eq:extended-potential}, any Hölder, Lipschitz or semiconcavity bound
on the functions $\varphi_\ell\circ t_a$, uniform in $a\in\X$ and $\ell\in[0,1]$, holds for $\bar f$
with the constant doubled. In Sections~\ref{sec:reg}--\ref{sec:entropy}, we establish uniform regularity
estimates in charts for these compositions, with smoothing and correction as needed, and transfer
these estimates to the potentials to bound the entropy of their restrictions to $\X$.

\subsection{Charts and distance functions}\label{sec:reg}

We fix finitely many charts covering $\X$ and record the Lipschitz and
one-sided second-difference bounds for $t_a:=d_M(\cdot,a)^q$, $a\in\X$,
used in the rest of this section.

\begin{lemma}[Distance bounds in charts]\label{lem:riemR}
Let $M$ be a complete, connected, smooth Riemannian manifold without boundary, of dimension
$d\ge1$, with geodesic distance $d_M$. Let $\X\subseteq M$ be compact with $\diam\X\le1$, equipped
with $\dX=d_M|_\X$, and let $q\ge1$. There are smooth charts $\psi_1,\dots,\psi_{j_0}:B_2\to M$ and
constants $\rho\in(0,1/8]$, $C_0\ge1$ and $C_t,K_t\ge0$ with the following properties.
Let $U:=\bigcup_j\psi_j(B_2)$ and $\delta_0:=\rho/(8C_0)$. For $a\in\X$ and each $j$, write
$w_{a,j}(z):=d_M(\psi_jz,a)$ and $t_{a,j}(z):=t_a(\psi_jz)=w_{a,j}(z)^q$ for the functions
$d_M(\cdot,a)$ and $t_a$ in the chart $\psi_j$.
\begin{enumerate}[label=\textup{(\roman*)},nosep]
\item \emph{Covering.} $\X\subseteq\bigcup_j\psi_j(B_1)$ and $\diam(\X\cup U)\le2$.
\item \emph{Lipschitz bounds.} For all $a\in\X$, all $j$, and all $z,z'\in B_2$,
$d_M(\psi_jz,\psi_jz')\le C_0\|z-z'\|$ and $|t_{a,j}(z)-t_{a,j}(z')|\le C_t\|z-z'\|$.
\item \emph{Second-difference bounds.} For all $a\in\X$, all $j$, and all $z,v$ with $z\pm v\in B_2$ and $\|v\|\le\delta_0$,
\begin{align*}
w_{a,j}(z)<2\rho &\quad\Longrightarrow\quad \Delta_v^2t_{a,j}(z)\ge-K_t\|v\|^2,\\
w_{a,j}(z)>\rho &\quad\Longrightarrow\quad \Delta_v^2t_{a,j}(z)\le K_t\|v\|^2.
\end{align*}
\end{enumerate}
\end{lemma}

Lemma~\ref{lem:riemR} is proved in Appendix~\ref{app:riem}. The bounds for $w_{a,j}$ follow from
Hessian comparison within the injectivity radius and from smooth upper supports
$d_M(a,a')+d_M(a',\cdot)$, with $a'$ on a minimizing geodesic from $a$, and those for $t_{a,j}$
follow since $q\ge1$. The charts and constants depend only on $q$ and the compact subset $\X$ of
$M$. The covering and Lipschitz bounds imply $N_{\mathcal X}(\eps)\le j_0(6C_0/\eps)^d$ for
$\eps\in(0,1]$. Since $\diam(\X\cup U)\le2$, the functions $t_a$ take values in $[0,2^q]$ on
$\X\cup U$, and $|t_a(x)-t_a(x')|\le q2^{q-1}d_M(x,x')$ for $x,x'\in\X\cup U$.

Henceforth, we fix charts and constants as in Lemma~\ref{lem:riemR}, with $U$, $\delta_0$, $w_{a,j}$ and
$t_{a,j}$ as defined there.

\subsection{\texorpdfstring{The case $p\le1$}{The case p <= 1}}\label{sec:pleq1}

This subsection proves Proposition~\ref{prop:Xside} for $p\le1$.
We use the Lipschitz bounds of Lemma~\ref{lem:riemR} to establish uniform
$p$-Hölder bounds for the potentials in charts, and then apply the entropy
bound for Hölder classes.

\begin{lemma}[Entropy of Hölder potentials]\label{lem:entHolder}
Let $p\le1$. For $\eps\in(0,1]$,
$\log N(\eps,\mathcal F_\X,\|\cdot\|_{\infty,\X})\le C\eps^{-d/p}$.
\end{lemma}

\begin{proof}
We follow the chartwise covering argument in the proof of
Lemma~A.4 of~\cite{hundrieser2024empirical}.
For $x,x'\in\X\cup U$, $a\in\X$ and $\ell\in[0,1]$,
$|\varphi_\ell(t_a(x))-\varphi_\ell(t_a(x'))|\le|t_a(x)-t_a(x')|^p\le(q2^{q-1}d_M(x,x'))^p$.
Integrating in~\eqref{eq:extended-cost} and then taking the minimum in~\eqref{eq:extended-potential}
gives the $Cd_M^p$-Hölder bound for each $\bar f$ on $\X\cup U$. Since
$|\bar f|=|f|\le2$ on $\X$ and $\diam(\X\cup U)\le2$, each $\bar f$ is also bounded on $\X\cup U$ by a
constant, and by the chart Lipschitz bound each $\bar f\circ\psi_j$ is $C\|z-z'\|^p$-Hölder on
$\bar B_1$. The standard sup-norm entropy bound for bounded $p$-Hölder functions on a ball is
$C\eps^{-d/p}$ \citep[Theorem~2.7.1]{vaart1996weak}; taking the product over charts gives the same
bound for $\mathcal F_\X=\{\bar f|_\X:(\gamma,f,g)\in\mathcal T\}$.
\end{proof}

\begin{proof}[Proof of Proposition~\ref{prop:Xside} for $p\le1$]
The class $\mathcal F_\X$ is uniformly bounded by~$2$ (Lemma~\ref{lem:pot}), so
Dudley's entropy bound (Lemma~\ref{lem:dudley}) with $a=p$ and Lemma~\ref{lem:entHolder} give
\[
  \E\sup_{f\in\mathcal F_\X}\Bigl|\int f\,d(\muh-\mu)\Bigr|\le C\psi_{d,p}(n)=C\psi_{d,p\wedge2}(n).\qedhere
\]
\end{proof}

\subsection{\texorpdfstring{The case $p>1$: smoothing and correction}{The case p > 1: smoothing and correction}}\label{sec:recentred}

Throughout this subsection and the next, $p>1$. Our goal is to use
semiconcavity in place of the Hölder continuity used in Section~\ref{sec:pleq1} to prove Proposition~\ref{prop:Xside}.
However, the compositions $\varphi_\ell\circ t_a$ need not be uniformly
semiconcave (see Remark~\ref{rem:correction}). We therefore smooth $\varphi_\ell$
when $1<p<2$ and, for all $p>1$, subtract from the cost a function of $x$ alone,
as outlined in Step~2 of Section~\ref{sec:overview}. Lemma~\ref{lem:key} shows
that the resulting compositions are semiconcave in charts, uniformly in $a$
and $\ell$.

\medskip\noindent\emph{Smoothing.}
We first replace $\varphi_\ell$ by a smooth approximation, which is needed for
$1<p<2$, where $\varphi_\ell$ is not semiconcave on any neighborhood of $\ell$. For $\ell\in[0,1]$ and $h\ge0$, let $\varphi_{\ell,h}(t):=((t-\ell)^2+h^2)^{p/2}$ for $t\ge0$, so that
$\varphi_{\ell,0}=\varphi_\ell$. Let $\mathsf H_p:=(0,1]$ if $1<p<2$ and $\mathsf H_p:=\{0\}$ if $p\ge2$,
and let $L_p:=p2^{q(p-1)}$. For $p\ge2$, we thus take $h=0$, and no smoothing is needed.
For $\gamma\in\PP(\X\times\Y)$ and $h\in\mathsf H_p$, and for $f,g$ with
$(\gamma,f,g)\in\mathcal T$, the smoothed cost and the smoothed potential are
\begin{equation}\label{eq:smoothed}
  c_{\gamma,h}(x,y):=2\int\varphi_{\dY(y,b)^q,h}\bigl(t_a(x)\bigr)\,d\gamma(a,b),\qquad
  f_h(x):=\min_{y\in\Y}\{c_{\gamma,h}(x,y)-g(y)\},
\end{equation}
for $x\in M$ and $y\in\Y$. We now bound the derivatives of $\varphi_{\ell,h}$ and show that
smoothing changes the potentials by at most $2h^p$.

\begin{lemma}[Smoothing]\label{lem:smooth}
Let $\ell\in[0,1]$ and $h\in\mathsf H_p$, and let $\lambda_h:=p(p-1)2^{q(p-2)}$ if $p\ge2$ and
$\lambda_h:=ph^{p-2}$ if $1<p<2$.
\begin{enumerate}[label=\textup{(\alph*)},nosep]
\item For $t\in[0,2^q]$,
$|\varphi'_{\ell,h}(t)|\le L_p$, $\varphi''_{\ell,h}(t)\le\lambda_h$ and $0\le\varphi_{\ell,h}(t)-\varphi_\ell(t)\le h^p$.
\item For $(\gamma,f,g)\in\mathcal T$, $0\le f_h-\bar f\le2h^p$ on $\X\cup U$.
\end{enumerate}
\end{lemma}

\begin{proof}
(a) Write $z:=t-\ell$, so that $|z|\le2^q$. Then
$\varphi'_{\ell,h}(t)=pz(z^2+h^2)^{p/2-1}$ and
$\varphi''_{\ell,h}(t)=p(z^2+h^2)^{p/2-2}\bigl((p-1)z^2+h^2\bigr)$. If $p\ge2$, then $h=0$, and both
derivative bounds follow from $|z|\le2^q$. If $1<p<2$, then $|\varphi'_{\ell,h}(t)|\le p|z|^{p-1}\le L_p$
and $\varphi''_{\ell,h}(t)\le p(z^2+h^2)^{p/2-1}\le ph^{p-2}$. The last bound in (a) follows from the
subadditivity of $x\mapsto x^{p/2}$ if $1<p<2$, and is trivial if $h=0$.
(b) Since $t_a\in[0,2^q]$ on $\X\cup U$, applying part (a) to the integrals
in~\eqref{eq:extended-cost} and~\eqref{eq:smoothed} gives $0\le c_{\gamma,h}-c_\gamma\le2h^p$ on
$(\X\cup U)\times\Y$, and taking the minima defining $\bar f$ and $f_h$
in~\eqref{eq:extended-potential} and~\eqref{eq:smoothed} gives $0\le f_h-\bar f\le2h^p$.
\end{proof}

\medskip\noindent\emph{Correction.}
We next construct a correction function $\Omega:\R\to\R$ so that the compositions $(\varphi_{\ell,h}-\Omega)\circ t_a$ are
semiconcave in charts, uniformly in $a$ and $\ell$. To this end, we choose $\Omega$ so that
$\varphi_{\ell,h}-\Omega$ is nonincreasing for small $t$ and nondecreasing for larger $t$, for all
$\ell$ and $h$. This allows us to use the lower second-difference bound for $t_a$ near $a$ and the
upper bound away from $a$ in Lemma~\ref{lem:riemR}. Let $\tau_1:=(9\rho/8)^q$ and $\tau_2:=(15\rho/8)^q$. Thus $t_a(x)\in[\tau_1,\tau_2]$ if and only if
$9\rho/8\le d_M(x,a)\le15\rho/8$, a range inside $(\rho,2\rho)$, where both second-difference bounds of
Lemma~\ref{lem:riemR} hold. Fix a nonincreasing $\chi\in C^1(\R;[-1,1])$, depending only on $\rho$
and $q$, with $\chi=1$ on $(-\infty,\tau_1]$ and $\chi=-1$ on $[\tau_2,\infty)$, and
define the correction function
\[
  \Omega(t):=L_p\int_0^t\chi(u)\,du,\qquad C_\Omega:=L_p\sup|\chi'| ,
\]
so that $\Omega\in C^2(\R)$, $\Omega'=L_p\chi$ and $|\Omega''|\le C_\Omega$. By
Lemma~\ref{lem:smooth}(a), each $\varphi_{\ell,h}-\Omega$ is nonincreasing on $[0,\tau_1]$ and
nondecreasing on $[\tau_2,2^q]$. For $\gamma\in\PP(\X\times\Y)$, define the corresponding correction by
\[
  \Omega_\gamma(x):=2\int\Omega\bigl(t_a(x)\bigr)\,d\gamma(a,b),\qquad x\in M.
\]
For $h\in\mathsf H_p$, and for $f,g$ with
$(\gamma,f,g)\in\mathcal T$, we call
\begin{equation}\label{eq:corrected}
  c_\gamma^\circ:=c_\gamma-\Omega_\gamma,\qquad c_{\gamma,h}^\circ:=c_{\gamma,h}-\Omega_\gamma,\qquad
  f^\circ:=\bar f-\Omega_\gamma,\qquad f_h^\circ:=f_h-\Omega_\gamma
\end{equation}
the corrected cost, the smoothed corrected cost, the corrected potential and the smoothed corrected
potential, respectively.
The costs are defined on $M\times\Y$, and the potentials on $M$.

The following elementary inequality will allow us to combine the monotonicity of
$\varphi_{\ell,h}-\Omega$ with the one-sided second-difference bounds for $t_a$ in
Lemma~\ref{lem:riemR} to prove semiconcavity of the compositions $(\varphi_{\ell,h}-\Omega)\circ t_a$ in
charts.

\begin{lemma}[Second-difference bound for compositions]\label{lem:comp}
Let $I$ be an interval and $G\in C^2(I)$ with $G''\le\lambda$ and $|G'|\le L_\ast$ on $I$. Let
$t_-,t_0,t_+\in I$ and $\Delta:=t_++t_--2t_0$. Then
\[
  G(t_+)+G(t_-)-2G(t_0)\le\tfrac\lambda4(t_+-t_-)^2+L_\ast\,\Delta_\sharp,
\]
where $\Delta_\sharp:=|\Delta|$ in general, $\Delta_\sharp:=\Delta_-$ if $G'\le0$ on $I$, and
$\Delta_\sharp:=\Delta_+$ if $G'\ge0$ on $I$; here $\Delta_\pm:=\max\{\pm\Delta,0\}$.
\end{lemma}

\begin{proof}
Let $\bar t:=(t_++t_-)/2\in I$. Since $G-\lambda t^2/2$ is concave,
$G(t_+)+G(t_-)-2G(\bar t)\le\lambda(t_+-t_-)^2/4$. Also $\bar t-t_0=\Delta/2$, so
$2(G(\bar t)-G(t_0))\le L_\ast|\Delta|$. If $G'\le0$, this difference is $\le0$ when $\Delta\ge0$, so it
is at most $L_\ast\Delta_-$; if $G'\ge0$, it is $\le0$ when $\Delta\le0$, so it is at most $L_\ast\Delta_+$.
Adding the bound from concavity to the bound on $2(G(\bar t)-G(t_0))$ gives the stated inequality.
\end{proof}

We now show that the smoothed and corrected compositions are uniformly Lipschitz and semiconcave
in charts.

\begin{lemma}[Semiconcavity after correction]\label{lem:key}
For $a\in\X$, $\ell\in[0,1]$ and $h\in\mathsf H_p$, let
\[
  u_{a,\ell,h}:=(\varphi_{\ell,h}-\Omega)\circ t_a\quad\text{on }\X\cup U .
\]
Then, for every $j$, the function $u_{a,\ell,h}\circ\psi_j$ is $2L_pC_t$-Lipschitz and
$\Lambda_h$-semiconcave on $B_2$, with
\[
  \Lambda_h:=(\lambda_h+C_\Omega)\,C_t^2+2L_pK_t .
\]
Thus $\Lambda_h\le C$ if $p\ge2$, and $\Lambda_h\le C(1+h^{p-2})$ if $1<p<2$.
\end{lemma}

\begin{proof}
Let $G:=\varphi_{\ell,h}-\Omega$. By Lemma~\ref{lem:smooth}(a), $|G'|\le2L_p$ and
$G''\le\lambda_h+C_\Omega$ on $[0,2^q]$.

\medskip\noindent\emph{Lipschitz bound.} By Lemma~\ref{lem:riemR}, $t_{a,j}$ is $C_t$-Lipschitz on
$B_2$. Together with $|G'|\le 2L_p$, this shows that
$u_{a,\ell,h}\circ\psi_j=G\circ t_{a,j}$ is $2L_pC_t$-Lipschitz on $B_2$.

\medskip\noindent\emph{Semiconcavity.} The correction function gives $G'=\varphi_{\ell,h}'-L_p\le0$ on
$[0,\tau_1]$ and $G'=\varphi_{\ell,h}'+L_p\ge0$ on $[\tau_2,2^q]$. To bound the second differences, fix $z\pm v\in B_2$ with $\|v\|\le\delta_0$, and write
$t_\pm:=t_{a,j}(z\pm v)$, $t_0:=t_{a,j}(z)$, $\Delta_v^2t:=t_++t_--2t_0$, $w_\pm:=w_{a,j}(z\pm v)$
and $w_0:=w_{a,j}(z)$. By the Lipschitz bounds, $|t_+-t_-|\le2C_t\|v\|$ and
$|w_\pm-w_0|\le C_0\delta_0=\rho/8$. We consider three cases according to $w_0=d_M(\psi_jz,a)$. In
cases (i) and (iii), the bound $|w_\pm-w_0|\le\rho/8$ determines the intervals containing
$t_-,t_0,t_+$.
\begin{enumerate}[label=(\roman*),nosep]
\item $w_0\le\rho$. The values $t_-,t_0,t_+$ all lie in $[0,\tau_1]$, where $G'\le0$, and $\Delta_v^2t\ge-K_t\|v\|^2$.
\item $\rho<w_0<2\rho$. Then $|\Delta_v^2t|\le K_t\|v\|^2$.
\item $w_0\ge2\rho$. The values $t_-,t_0,t_+$ all lie in $[\tau_2,2^q]$, where $G'\ge0$, and $\Delta_v^2t\le K_t\|v\|^2$.
\end{enumerate}
In case (i), the lower bound combines with $G'\le0$, and in case (iii), the upper bound combines
with $G'\ge0$. In case (ii), both bounds hold, so no sign condition on $G'$ is needed.
In each case, Lemma~\ref{lem:comp}, with $\Delta_\sharp\le K_t\|v\|^2$ and
$(t_+-t_-)^2\le4C_t^2\|v\|^2$, gives $\Delta_v^2(u_{a,\ell,h}\circ\psi_j)(z)\le\Lambda_h\|v\|^2$.
Since this bound holds for all $\|v\|\le\delta_0$, the function $u_{a,\ell,h}\circ\psi_j$ is
$\Lambda_h$-semiconcave on $B_2$.
\end{proof}

\begin{remark}[Failure of semiconcavity without correction]\label{rem:correction}
The following example shows that an uncorrected potential and a composition
$\varphi_\ell\circ t_a$ can both fail to be semiconcave, even on a circle. Let $\X$ be a circle of circumference $2r<2$, $\Y=\{b,y\}$ with $\dY(b,y)=1$,
$\gamma:=\delta_{(a,b)}$ with $a\in\X$, and let $g_0(y)=0$, $g_0(b)<-2$. Then
$f:=\min_{y'\in\Y}(c_\gamma(\cdot,y')-g_0(y'))=2\varphi_1\circ t_a=2(1-t_a)^p$, and $(f,g)\in\Phi_{\X,\Y}(c_\gamma)$ for
$g(y'):=\min_{x\in\X}(c_\gamma(x,y')-f(x))$.

Near the antipodal point of $a$, using signed distance $z$, we have
$t_a=(r-|z|)^q$, so
$\Delta_v^2t_a(0)=-2qr^{q-1}|v|+O(v^2)$ as $v\to0$.
Since $\varphi_1$ is smooth near $r^q<1$ and $\varphi_1'(r^q)<0$,
\[
\begin{aligned}
  \Delta_v^2f(0)
  &=2\varphi_1'(r^q)\Delta_v^2t_a(0)+O(v^2)\\
  &=4pq\,r^{q-1}(1-r^q)^{p-1}|v|+O(v^2).
\end{aligned}
\]
The positive linear term rules out semiconcavity. This example also illustrates the second obstacle mentioned in Step~2 of Section~\ref{sec:overview}.
\end{remark}

\subsection{\texorpdfstring{The case $p>1$: entropy and rate}{The case p > 1: entropy and rate}}\label{sec:entropy}

Using the semiconcavity established in Lemma~\ref{lem:key}, we apply Bronshtein's entropy bound
(Theorem~\ref{thm:bron}) to the smoothed corrected potentials and control the smoothing error to
obtain an entropy bound for the corrected potentials. Dudley's entropy bound and a separate estimate for the corrections then complete
the proof of Proposition~\ref{prop:Xside} for $p>1$.

Let $\mathcal E$ be the class of corrected potentials on $\X$, and write $\mathcal W$ for the family of
the functions $\Omega\circ t_a$, whose mixtures give the corrections $\Omega_\gamma$:
\[
  \mathcal E:=\bigl\{f^\circ|_\X:(\gamma,f,g)\in\mathcal T\bigr\},\qquad
  \mathcal W:=\bigl\{\Omega\circ t_a|_\X:a\in\X\bigr\}.
\]
Entropy bounds for semiconcave functions reduce to the following bound for convex functions.

\begin{theorem}[Bronshtein's entropy bound; {\citep{bronshtein1976entropy}; \citep[Corollary~2.7.10]{vaart1996weak}; see also \citep[Theorem~3.2]{guntuboyina2013covering}}]\label{thm:bron}
Let $d\ge1$, let $Z\subseteq\R^d$ be compact and convex, and let $\mathcal C$ be the class of convex
$L_\ast$-Lipschitz functions $F:Z\to[0,1]$. Then
$\log N(\eps,\mathcal C,\|\cdot\|_\infty)\le C_{d,Z}(1+L_\ast)^{d/2}\eps^{-d/2}$.
\end{theorem}

Combined with Lemma~\ref{lem:key}, it gives the entropy bound for $\mathcal E$ in the following lemma.
The bound for $\mathcal W$ does not use Lemma~\ref{lem:key} and only relies on the Lipschitz dependence
of $\Omega\circ t_a$ on $a$.

\begin{lemma}[Entropy of corrected potentials]\label{lem:entE}
For $\eps\in(0,1]$,
\[
  \log N(\eps,\mathcal E,\|\cdot\|_{\infty,\X})\le C\eps^{-d/(p\wedge2)},\qquad
  \log N(\eps,\mathcal W,\|\cdot\|_{\infty,\X})\le d\log(C/\eps).
\]
\end{lemma}

\begin{proof}
\emph{Entropy of $\mathcal E$.} Set $L_\ast:=4L_pC_t$ and $B_0:=2(4^q+1)^{p/2}+4+2^{q+1}L_p$. For $\lambda\ge0$ let
$\mathcal S(\lambda)$ be the class of $e:\X\cup U\to[-B_0,B_0]$ for which each $e\circ\psi_j$ is
$L_\ast$-Lipschitz and $\lambda$-semiconcave on $B_2$. Applying Bronshtein's bound in each chart, as in
the proof of~\citep[Lemma~A.3]{hundrieser2024empirical}, gives
\begin{equation}\label{eq:semiconcave-entropy}
  \log N\bigl(2\eps,\mathcal S(\lambda)|_\X,\|\cdot\|_{\infty,\X}\bigr)
  \le Cj_0(1+\lambda)^{d/2}\eps^{-d/2}.
\end{equation}
Indeed, for $e\in\mathcal S(\lambda)$, the function $F:=\lambda\|\cdot\|^2/2-e\circ\psi_j+B_0$ is convex
on $\bar B_1$ with values in $[0,b_\lambda]$, where $b_\lambda:=\lambda/2+2B_0$, and $F/b_\lambda$ is
Lipschitz with a constant independent of $\lambda$. Applying Theorem~\ref{thm:bron} to $F/b_\lambda$
at radius $\eps/b_\lambda$ in each chart gives the factor $b_\lambda^{d/2}\le C(1+\lambda)^{d/2}$. The
factor $2$ in~\eqref{eq:semiconcave-entropy} comes from combining the charts and passing to an
internal cover.

Now let $(\gamma,f,g)\in\mathcal T$ and $h\in\mathsf H_p$. Since $\Omega_\gamma$ does not depend on $y$,
the smoothed corrected potential and cost satisfy
\[
  f_h^\circ(x)=\min_{y\in\Y}\{c_{\gamma,h}^\circ(x,y)-g(y)\},\qquad
  c_{\gamma,h}^\circ(x,y)=2\int u_{a,\dY(y,b)^q,h}(x)\,d\gamma(a,b),
\]
so $f_h^\circ$ is a pointwise minimum of mixtures of the functions $2u_{a,\ell,h}$ minus constants.
By the reduction to regularity in Section~\ref{sec:notation}, Lemma~\ref{lem:key} therefore shows that the composition of
$f_h^\circ$ with each chart $\psi_j$ is $2\Lambda_h$-semiconcave and $L_\ast$-Lipschitz on $B_2$. Moreover, $0\le c_{\gamma,h}\le2(4^q+1)^{p/2}$ and
$|\Omega_\gamma|\le2^{q+1}L_p$ on $\X\cup U$, and $|g|\le4$. Thus $f_h^\circ|_{\X\cup U}\in\mathcal S(2\Lambda_h)$, with
the envelope $B_0$ and Lipschitz constant $L_\ast$ independent of $h$ and $\gamma$. By Lemma~\ref{lem:smooth}(b),
$\|f_h^\circ-f^\circ\|_{\infty,\X\cup U}\le2h^p$.

Hence a cover of $\mathcal S(2\Lambda_h)|_\X$ yields a cover of the corrected potentials, once $h$
is chosen so that the smoothing error $2h^p$ is at most half the desired covering radius.
If $p\ge2$, take $h=0$. Then $\mathcal E\subseteq\mathcal S(2\Lambda_0)|_\X$, so~\eqref{eq:semiconcave-entropy}
gives the entropy bound $C\eps^{-d/2}$. If $1<p<2$, take $h:=(\eps/4)^{1/p}\in(0,1]$. Then $2h^p=\eps/2$ and
$1+2\Lambda_h\le C\eps^{(p-2)/p}$. The bound~\eqref{eq:semiconcave-entropy}, applied with $\eps/4$ in place of $\eps$, gives an $\eps/2$-cover of
$\mathcal S(2\Lambda_h)|_\X$, which is an external $\eps$-cover of $\mathcal E$, of log-size at most
$C(\eps^{(p-2)/p}\eps^{-1})^{d/2}=C\eps^{-d/p}$. Passing to an internal cover changes only $C$.

\medskip\noindent\emph{Entropy of $\mathcal W$.} The map $a\mapsto\Omega\circ t_a$ is $L_p\,q\,2^{q-1}$-Lipschitz
in the sup norm. An $\eps$-net of $\X$ therefore induces a sup-norm net of $\mathcal W$, with its
radius multiplied by this Lipschitz constant. The bound $N_{\mathcal X}(\eps)\le j_0(6C_0/\eps)^d$
gives $\log N(\eps,\mathcal W,\|\cdot\|_{\infty,\X})\le d\log(C/\eps)$.
\end{proof}

\begin{proof}[Proof of Proposition~\ref{prop:Xside} for $p>1$]
On $\X$, each $f\in\mathcal F_\X$ equals $f^\circ+\Omega_\gamma$ with $f^\circ|_\X\in\mathcal E$, and by Fubini,
\[
  \sup_\gamma\Bigl|\int \Omega_\gamma\,d(\muh-\mu)\Bigr|\le2\sup_{w\in\mathcal W}\Bigl|\int w\,d(\muh-\mu)\Bigr| .
\]
The classes $\mathcal E$ and $\mathcal W$ are uniformly bounded, so applying Dudley's entropy bound
(Lemma~\ref{lem:dudley}) to both, with the entropy bounds of Lemma~\ref{lem:entE}, gives
\[
  \E\sup_{f\in\mathcal F_\X}\Bigl|\int f\,d(\muh-\mu)\Bigr|
  \le\underbrace{\E\sup_{u\in\mathcal E}\Bigl|\int_\X u\,d(\muh-\mu)\Bigr|}_{\le\,C\psi_{d,p\wedge2}(n)}
  +\underbrace{2\,\E\sup_{w\in\mathcal W}\Bigl|\int w\,d(\muh-\mu)\Bigr|}_{\le\,Cn^{-1/2}}
  \le C\psi_{d,p\wedge2}(n),
\]
where the last step uses $n^{-1/2}\lesssim\psi_{d,p\wedge2}(n)$.
\end{proof}

\subsection{\texorpdfstring{Potentials on $\Y$}{Potentials on Y}}\label{sec:yside}

This subsection proves Proposition~\ref{prop:Yside}. The argument on $\X$ used regularity in
charts; here we express each $g\in\mathcal F_\Y$ as the $c_\gamma^\circ$-transform of a corrected
potential on $\X$; see~\eqref{eq:gtransform}. Discretizing the minimum defining this transform allows us to construct a cover
of $\mathcal F_\Y$ from covers of the corrected potentials and the costs restricted to the net
(Lemma~\ref{lem:discrete}).
We thereby bound its entropy in $L^2(Q)$, uniformly over finitely supported probability measures
$Q$ on $\Y$.

To treat the cases $p\le1$ and $p>1$ together, set $\Omega:=0$ if $p\le1$, so that $\Omega_\gamma=0$,
$c_\gamma^\circ=c_\gamma$, $f^\circ=\bar f$ and $\mathcal E=\mathcal F_\X$ in that case. By Lemmas~\ref{lem:entHolder} and~\ref{lem:entE}, recall
\begin{equation}\label{eq:Eentropy}
  \log N(\eps,\mathcal E,\|\cdot\|_{\infty,\X})\le C\eps^{-d/(p\wedge2)}
\end{equation}
for all
$p>0$ and $\eps\in(0,1]$.
For $(\gamma,f,g)\in\mathcal T$, we have $f^\circ|_\X=f-\Omega_\gamma\in\mathcal E$, and since
$\Omega_\gamma$ does not depend on $y$, the conjugacy relations give
\begin{equation}\label{eq:gtransform}
  g(y)=\min_{x\in\X}\{c_\gamma^\circ(x,y)-f^\circ(x)\},\qquad y\in\Y.
\end{equation}

We first collect the bounds that control the error from discretizing the minimum
in~\eqref{eq:gtransform} and the entropy of the discretized costs.
Write $\theta=(a,b)\in\Theta:=\X\times\Y$, and let $c_\theta^\circ:=c_{\delta_\theta}^\circ$
be the corrected cost for the point mass at $\theta$, so that, for $x\in\X$ and
$y\in\Y$,
\[
  c_\theta^\circ(x,y)=2\bigl[\varphi_{\dY(y,b)^q}(t_a(x))-\Omega(t_a(x))\bigr],\qquad
  c_\gamma^\circ(x,y)=\int c_\theta^\circ(x,y)\,d\gamma(\theta).
\]
By Lemmas~\ref{lem:reg} and~\ref{lem:pot} and since $|\Omega'|\le L_p$ when $p>1$,
there is $L_c\ge1$ such that, for $\theta=(a,b)$ and $\theta'=(a',b')$ in $\Theta$,
\begin{enumerate}[label=\textup{(\alph*)},nosep]
\item $|c_\theta^\circ|\le L_c$ on $\X\times\Y$;
\item $|c_\theta^\circ-c_{\theta'}^\circ|\le L_c(d_M(a,a')+\dY(b,b'))^{p\wedge1}$ on $\X\times\Y$;
\item every $x\mapsto c_\gamma^\circ(x,y)$ and every $e\in\mathcal E$ is
      $L_cd_M^{p\wedge1}$-Hölder on $\X$.
\end{enumerate}

\medskip\noindent\emph{Reduction to the corrected costs.}
We replace the minimum over $\X$ in~\eqref{eq:gtransform} by a minimum over a finite net, which
reduces covering $\mathcal F_\Y$ to covering $\mathcal E$ and the corrected costs restricted to this net in their first argument.

In the rest of this subsection, $Q$ is a finitely supported probability measure on $\Y$ and
$\eps\in(0,1]$. Set
$\eta:=(\eps/(6L_c))^{1/(p\wedge1)}$ and fix an $\eta$-net $\{x_1,\ldots,x_k\}\subseteq\X$ with
$k\le j_0(6C_0/\eta)^d\le C\eps^{-d/(p\wedge1)}$. For $\gamma\in\PP(\X\times\Y)$, let
$H_\gamma(y,i):=c_\gamma^\circ(x_i,y)$ be the discretized corrected cost on this net.
For functions $H$ on $\Y\times\{1,\ldots,k\}$ and $e$ on $\X$, let
\[
  G_{H,e}(y):=\min_{1\le i\le k}\{H(y,i)-e(x_i)\},\qquad
  \|H\|_{Q,\infty}:=\Bigl(\int_\Y\max_{1\le i\le k}|H(y,i)|^2\,dQ(y)\Bigr)^{1/2}.
\]
We now show
that discretizing the minimum over $\X$ costs at most $\eps/3$ and that $G_{H,e}$ depends Lipschitz
continuously on $(H,e)$.

\begin{lemma}[Discretization and stability of $c$-transforms]\label{lem:discrete}
\begin{enumerate}[label=\textup{(\roman*)},nosep]
\item For $(\gamma,f,g)\in\mathcal T$, we have $0\le G_{H_\gamma,f^\circ}-g\le\eps/3$
on $\Y$.
\item For all $H,H'$ and $e,e'$ as above,
$\|G_{H,e}-G_{H',e'}\|_{L^2(Q)}\le\|e-e'\|_{\infty,\X}+\|H-H'\|_{Q,\infty}$.
\end{enumerate}
\end{lemma}

\begin{proof}
(i) Since the net lies in $\X$,~\eqref{eq:gtransform} gives $G_{H_\gamma,f^\circ}\ge g$. Conversely, for
$y\in\Y$ let $x\in\X$ attain the minimum in~\eqref{eq:gtransform}, and choose $x_i$ with
$d_M(x,x_i)\le\eta$. Property~\textup{(c)} gives
$G_{H_\gamma,f^\circ}(y)\le c_\gamma^\circ(x_i,y)-f^\circ(x_i)\le g(y)+2L_c\eta^{p\wedge1}=g(y)+\eps/3$.
(ii) The elementary inequality $|\min_iu_i-\min_iv_i|\le\max_i|u_i-v_i|$ gives
$|G_{H,e}(y)-G_{H',e'}(y)|\le\|e-e'\|_{\infty,\X}+\max_i|H(y,i)-H'(y,i)|$, and taking the
$L^2(Q)$ norm gives the claim.
\end{proof}

Part~\textup{(i)} bounds the error from discretizing the $c$-transform
by $\eps/3$. Part~\textup{(ii)} then shows that $\eps/3$-covers of
$\mathcal E$ in $\|\cdot\|_{\infty,\X}$ and of $\{H_\gamma\}$ in
$\|\cdot\|_{Q,\infty}$ yield an external $\eps$-cover of $\mathcal F_\Y$
in $L^2(Q)$. Since the entropy of $\mathcal E$ is already controlled
by~\eqref{eq:Eentropy}, it remains to bound the entropy of $\{H_\gamma\}$.

\medskip\noindent\emph{Entropy of the corrected costs.}
For $\theta\in\Theta$, let $H_\theta(y,i):=c_\theta^\circ(x_i,y)$. Let $\mathrm{Unif}_k$ be the uniform probability measure
on $\{1,\ldots,k\}$, and let $P:=Q\otimes \mathrm{Unif}_k$. For $r\ge2$, the comparison of
finite-dimensional $\ell^r$ and $\ell^\infty$ norms followed by Jensen's inequality gives
\begin{equation}\label{eq:cost-norm-comparison}
  \|H\|_{Q,\infty}\le k^{1/r}\|H\|_{L^r(P)}.
\end{equation}
This comparison allows us to use an $L^r(P)$ entropy bound for the discretized costs.
Since each $H_\gamma$ is a mixture of the maps $H_\theta$, the family
$\{H_\gamma\}$ lies in the closed convex hull of $\{H_\theta:\theta\in\Theta\}$
in $L^r(P)$. We thus use the following entropy bound for convex hulls.

\begin{theorem}[$L^r$ convex-hull entropy; {\citep{carl1999metric}; \citep[Theorem~2.1, eq.~(2.3)]{song2002upper}}]\label{thm:carl}
Fix $r\in[2,\infty)$, $D>0$ and $C_{\mathcal A}\ge1$. Let $P$ be any probability
measure and let $\mathcal A$ be a class of measurable real-valued functions
with $|F|\le B$ for every $F\in\mathcal A$, where $B>0$. If
\[
  N(\delta B,\mathcal A,L^r(P))\le(C_{\mathcal A}/\delta)^D
  \qquad (0<\delta\le1),
\]
then, writing $\kappa:=2D/(D+2)<2$,
\[
  \log N(\delta B,\operatorname{conv}\mathcal A,L^r(P))
  \le C\delta^{-\kappa}\qquad(0<\delta\le1),
\]
where $C$ depends only on $r,D,C_{\mathcal A}$. In particular, $C$ can be chosen independently of $P$ whenever
the assumed covering bound holds uniformly in $P$ with the same
$D$ and $C_{\mathcal A}$.
\end{theorem}

The entropy exponent $\kappa=2D/(D+2)$ in Theorem~\ref{thm:carl} is
independent of $r$. We therefore choose $r$ large enough
that the factor $k^{1/r}$ leaves the resulting exponent below $2$.

\begin{lemma}[Entropy of discretized costs]\label{lem:costs}
Let $D:=(d+d')/(p\wedge1)$ and $\kappa:=2D/(D+2)$, and fix $r\ge2$ with
$r>dD/(2(p\wedge1))$, so that
\[
  \kappa_r:=\kappa\left(1+\frac{d}{(p\wedge1) r}\right)<2.
\]
Then, with $C$ independent of $Q$ and $\eps$,
\[
  \log N\bigl(\eps/3,\{H_\gamma:\gamma\in\PP(\X\times\Y)\},\|\cdot\|_{Q,\infty}\bigr)\le C\eps^{-\kappa_r}.
\]
\end{lemma}

\begin{proof}
We have $H_\gamma=\int H_\theta\,d\gamma(\theta)$, where the integral is understood as the Bochner
integral of $\theta\mapsto H_\theta$ over $\gamma$. Hence $\{H_\gamma\}$ is contained in the closure of
$\operatorname{conv}\mathcal A_\eps$, where $\mathcal A_\eps:=\{H_\theta:\theta\in\Theta\}$, by the same
approximation argument as in~\citep[Appendix~D.5, Step~1]{leung2026sample}.
Since $P$ has finite support, the ambient $L^r(P)$ space is finite-dimensional. A finite cover of
$\operatorname{conv}\mathcal A_\eps$ by closed balls also covers its closure, and passing to an
internal cover at most doubles the radius.
Moreover, properties~\textup{(a)} and~\textup{(b)} give, uniformly in $Q$, $k$ and $\eps$,
\[
  \|H_\theta\|_{L^r(P)}\le L_c,\qquad
  \|H_\theta-H_{\theta'}\|_{L^r(P)}
  \le L_c(d_M(a,a')+d_\Y(b,b'))^{p\wedge1}.
\]
By the covering bound for $\X$ and the covering assumption on $\Y$, for $0<\delta\le1$,
\[
  N(\delta L_c,\mathcal A_\eps,L^r(P))
  \le N_\X(\delta^{1/(p\wedge1)}/2)
       N_\Y(\delta^{1/(p\wedge1)}/2)
  \le(C/\delta)^D.
\]
Theorem~\ref{thm:carl}, with the fixed $r$ above, therefore gives
\[
  \log N\bigl(t,\{H_\gamma\},L^r(P)\bigr)
  \le C_r(L_c/t)^\kappa \qquad (0<t\le L_c).
\]
By the norm comparison~\eqref{eq:cost-norm-comparison}, an $\eps/(3k^{1/r})$-cover in $L^r(P)$ is an $\eps/3$-cover in $\|\cdot\|_{Q,\infty}$, and
\[
\begin{aligned}
  \log N\bigl(\eps/3,\{H_\gamma\},\|\cdot\|_{Q,\infty}\bigr)
  &\le C_r\left(\frac{3L_c k^{1/r}}{\eps}\right)^\kappa\\
  &\le C\eps^{-\kappa(1+d/((p\wedge1) r))}
   =C\eps^{-\kappa_r}.
\end{aligned}
\]
The condition $r>dD/(2(p\wedge1))$ is equivalent to $\kappa_r<2$.
\end{proof}

\medskip\noindent\emph{Entropy and rate.}
We now combine Lemmas~\ref{lem:discrete} and~\ref{lem:costs} with the entropy
bound~\eqref{eq:Eentropy} for the corrected potentials.

\begin{lemma}[Entropy of potentials on $\Y$]\label{lem:entY}
With $\kappa_r<2$ as in Lemma~\ref{lem:costs}, for every finitely supported probability measure
$Q$ on $\Y$ and every $\eps\in(0,1]$,
\[
  \log N\bigl(\eps,\mathcal F_\Y,L^2(Q)\bigr)
  \le C\bigl[\eps^{-\kappa_r}+\eps^{-d/(p\wedge2)}\bigr].
\]
\end{lemma}

\begin{proof}
Take an $\eps/3$-cover of $\{H_\gamma\}$ in $\|\cdot\|_{Q,\infty}$ and an $\eps/3$-cover of $\mathcal E$
in $\|\cdot\|_{\infty,\X}$. For $g\in\mathcal F_\Y$, with $(\gamma,f,g)\in\mathcal T$, choose centers
$\widetilde H$ and $\widetilde e$ with
$\|H_\gamma-\widetilde H\|_{Q,\infty}\le\eps/3$ and $\|f^\circ-\widetilde e\|_{\infty,\X}\le\eps/3$. By
Lemma~\ref{lem:discrete},
\[
\begin{aligned}
  \|g-G_{\widetilde H,\widetilde e}\|_{L^2(Q)}
  &\le\|g-G_{H_\gamma,f^\circ}\|_{L^2(Q)}
     +\|f^\circ-\widetilde e\|_{\infty,\X}
     +\|H_\gamma-\widetilde H\|_{Q,\infty}\\
  &\le\eps/3+\eps/3+\eps/3=\eps.
\end{aligned}
\]
Thus the functions $G_{\widetilde H,\widetilde e}$ form an external $\eps$-cover of $\mathcal F_\Y$ in
$L^2(Q)$, of logarithmic size at most $C[\eps^{-\kappa_r}+\eps^{-d/(p\wedge2)}]$ by
Lemma~\ref{lem:costs} and~\eqref{eq:Eentropy}. Passing to an internal cover at most doubles the radius,
and replacing $\eps$ by $\eps/2$ proves the claim.
\end{proof}

\begin{proof}[Proof of Proposition~\ref{prop:Yside}]
The class $\mathcal F_\Y$ is uniformly bounded by~$4$
(Lemma~\ref{lem:pot}), so applying Dudley's entropy bound
(Lemma~\ref{lem:dudley}) with the entropy bound of
Lemma~\ref{lem:entY} gives
\[
  \E\sup_{g\in\mathcal F_\Y}
  \Bigl|\int g\,d(\hat\nu_m-\nu)\Bigr|
  \le C\bigl[\psi_{d,p\wedge2}(m)+m^{-1/2}\bigr]
  \le C\psi_{d,p\wedge2}(m).
\]
Here the entropy term $\eps^{-d/(p\wedge2)}$ gives
$\psi_{d,p\wedge2}(m)$, while $\eps^{-\kappa_r}$ contributes
$O(m^{-1/2})$ because $\kappa_r<2$. The last step uses
$m^{-1/2}\lesssim\psi_{d,p\wedge2}(m)$.
\end{proof}

\subsection{Completion of the upper bound}\label{sec:proof}

We combine Propositions~\ref{prop:stab} and~\ref{prop:potentials} to
prove~\eqref{eq:upper}. Apply Proposition~\ref{prop:stab} with $(\mu_0,\nu_0)=(\mu,\nu)$ and
$(\mu_1,\nu_1)=(\muh,\hat\nu_m)$ to obtain
\[
  \bigl|D_{p,q}^p(\muh,\hat\nu_m)-D_{p,q}^p(\mu,\nu)\bigr|\le\sup_{f\in\mathcal F_\X}\Bigl|\int f\,d(\muh-\mu)\Bigr|+\sup_{g\in\mathcal F_\Y}\Bigl|\int g\,d(\hat\nu_m-\nu)\Bigr| .
\]
Taking expectations and applying Proposition~\ref{prop:potentials} proves~\eqref{eq:upper}
for the normalized spaces. The scaling at the beginning of Section~\ref{sec:reduction} gives the
result for the original spaces.\qed
\par

\begingroup
\newcommand{\iid}{\mathrel{\stackrel{\mathrm{i.i.d.}}{\sim}}}
\newcommand{\Var}{\operatorname{Var}}
\newcommand{\Unif}{\operatorname{Unif}}
\section{Proof of the lower bound}
\label{app:lower-bound}

In this section, we prove the lower bound~\eqref{eq:lowerintro}.
By isometric invariance, it suffices to construct the measures on $B_1^d(0)$ with $d:=d_\star$,
and then push them forward under the canonical embeddings into $\R^{d_x}$ and $\R^{d_y}$,
as in~\citep[Section~5.2]{leung2026sample}. These embeddings preserve $D_{p,q}^p$ for both the population and empirical measures.
Set $N:=n\wedge m$, so that it suffices to prove a lower bound of $cN^{-(p\wedge2)/(d\vee2(p\wedge2))}$. We construct measures $\mu,\, \nu \in \PP(B_1^d(0))$ that attain the desired lower bound. Let us proceed by splitting into two cases based on the dimension $d$. When $d<2(p\wedge2)$, the target bound \eqref{eq:lowerintro} reduces to $cN^{-1/2}$. In this case, we prove the following proposition.

\begin{proposition}[Case $d<2(p\wedge2)$]
\label{prop:lower-small-dim}
    There exist a constant $c>0$
    and an integer $n_0$ such that, for every $n,m\ge n_0$, there exist probability measures $\mu,\nu \in \PP(B_1^d(0))$ satisfying
    \[
     \E\left|
     D_{p,q}^p(\hat\mu_n,\hat\nu_m)-D_{p,q}^p(\mu,\nu)
     \right|
     \ge cN^{-1/2}.
    \]
\end{proposition}

\begin{proof}[Proof of Proposition~\ref{prop:lower-small-dim}]
Let $e_1$ be the first coordinate vector, and fix
$\theta \in (0,1) \setminus \{1/2\}$. Define the two measures $\mu:=(1-\theta)\delta_0+\theta\delta_{e_1}$ and $\nu:=\delta_0$, supported on $\{0,e_1\}$ and $\{0\}$, respectively,
where $\delta_0$ and $\delta_{e_1}$ are Dirac masses. Then $D_{p,q}^p(\mu,\nu) = 2\theta(1-\theta)$.
For $X_1,\ldots,X_n\iid \mu$, define $B_i:=\mathbf1_{\{X_i=e_1\}}$ and $\hat\theta_n:=\frac1n\sum_{i=1}^n B_i$. Then $B_1,\ldots,B_n$ are independent Bernoulli$(\theta)$ variables, and
$\hat\mu_n=(1-\hat\theta_n)\delta_0+\hat\theta_n\delta_{e_1}$.
Since $\hat\nu_m=\delta_0$ and $D_{p,q}^p(\hat\mu_n,\hat\nu_m) = 2\hat\theta_n(1-\hat\theta_n)$, a simple computation gives $$D_{p,q}^p(\hat \mu_n,\hat \nu_m) - D_{p,q}^p(\mu,\nu) = 2(1-2\theta)U_n-2U_n^2,$$ where $U_n:=\hat\theta_n-\theta = \frac{1}{n} \sum_{i=1}^n (B_i - \theta)$. Now since $n \hat \theta_n \sim \mathrm{Bin}(n,\theta)$,
it follows from~\citep[Theorem~1]{berend2013sharp} that, for all sufficiently
large $n$,
\[
    \E|U_n|
    =\frac1n\E|n\hat\theta_n-n\theta|
    \geq\sqrt{\frac{\theta(1-\theta)}{2n}}
    \geq\frac{\theta(1-\theta)}{\sqrt n}.
\]
Also, $\E U_n^2=\theta(1-\theta)/n$ by independence.
Therefore,
\begin{align*}
    \E\left|
        D_{p,q}^p(\hat\mu_n,\hat\nu_m)-D_{p,q}^p(\mu,\nu)
    \right|
    &\geq
    \frac{2|1-2\theta|\theta(1-\theta)}{\sqrt n}
    -\frac{2\theta(1-\theta)}n \geq
    \frac{|1-2\theta|\theta(1-\theta)}{\sqrt n},
\end{align*}
where the last inequality holds for
$n\geq4/(1-2\theta)^2$. Interchanging the two marginals then gives $\E\left|D_{p,q}^p(\hat\mu_n,\hat\nu_m)-D_{p,q}^p(\mu,\nu)\right| \geq \frac{|1-2\theta|\theta(1-\theta)}{\sqrt{m}}$. Hence,
\[
    \sup_{\mu,\nu\in\PP(B_1^d(0))} \E\left|D_{p,q}^p(\hat\mu_n,\hat\nu_m)-D_{p,q}^p(\mu,\nu)\right| \geq |1-2\theta|\theta(1-\theta)\max\{n^{-1/2},m^{-1/2}\} = cN^{-1/2}
\]
for some $c>0$.
\end{proof}

It remains to consider the case $d \geq 2(p\wedge2)$, in which we prove the following proposition.

\begin{proposition}[Case $d \geq 2(p\wedge2)$]
\label{prop:lower-large-dim}
    There exist a constant $c>0$ and an integer $n_0$ such that, for every $n,m\ge n_0$, there exist probability
    measures $\mu,\nu \in \PP(B_1^d(0))$ satisfying
    \[
     \E\left|
     D_{p,q}^p(\hat\mu_n,\hat\nu_m)-D_{p,q}^p(\mu,\nu)
     \right|
     \ge cN^{-(p\wedge2)/d}.
    \]
\end{proposition}

Our proof mirrors that of Zhang et al.~\citep{zhang2024gromov} and follows in several steps; the main idea is to bound $\E\bigl|D_{p,q}^p(\hat\mu_n,\hat\nu_m)-D_{p,q}^p(\mu,\nu)\bigr|$ from below using a transport cost on a cube. To this end, we first construct measures $\mu$ and $\nu$ using a combination of delta masses and a uniform distribution and expand upon the reason for such a construction below. For $0<p<2$, $D_{p,q}^p(\mu,\nu)=0$, and we lower bound $D_{p,q}^p(\hat\mu_n,\hat\nu_m)$ by a transport cost on the cube. For $p\geq2$, we lower bound the plug-in error using a term $I_{p,q}$ that depends purely on the marginals and an auxiliary functional $\mathcal R_{p,q}$, defined below. We compute the expected marginal term and bound $\mathcal R_{p,q}$ below by a transport cost.

\paragraph{Step 1: constructing the measures.}

Let $e_1,\dots,e_d \in \R^d$ be the standard coordinate vectors in $\R^d$. Let $r_0:=0$ and $r_k := \frac{1}{2} + \frac{k}{4d}$ for $k = 1,\dots,d$, and let $a_0 := 0$, $a_k := r_ke_k$, and $A := \{a_0,a_1,\dots,a_d\}$. Also define the cube $Q := \prod_{i=1}^d \left[\frac{1}{8\sqrt{d}},\frac{3}{16\sqrt{d}} \right]$ and $H := A \cup Q$. Write $H_k:=\{a_k\}$ for $0\leq k\leq d$ and $H_{d+1}:=Q$. The cube has side length $\ell:=1/(16\sqrt d)$. We then define the measures
\[
    \mu = \frac{2}{3} \delta_{a_0} + \frac{1}{6d} \sum_{k=1}^d \delta_{a_k} + \frac{1}{6} \lambda_Q, \qquad \nu = S_\# \mu, \qquad S(x) = \begin{cases}
        x, \quad &0<p<2 \\
        x/2, \quad &p\geq 2
    \end{cases},
\]
where $\lambda_Q$ is the uniform probability measure on the cube $Q$. Figure~\ref{fig:rho} shows $\mu$ for $d=2$ and $d=3$. We provide some insights as to why such a construction makes sense. We use transport costs to prove the lower bound, yet $D_{p,q}^p$ only involves pairwise distances and is invariant under isometries. To account for this invariance, Zhang et al.~\citep{zhang2024gromov} take an infimum over isometries and, for measures with nondegenerate covariance, relate the $(2,2)$-GW distance to $\inf_{U\in E(d)}W_2(\mu,U_\#\nu)$ through a covariance-dependent factor, where $E(d)$ is the isometry group on $\R^d$. In our proof, we take a different route. The delta masses at $a_0,\dots,a_d$ constrain the couplings in our construction and allow us to bound $D_{p,q}^p$ for $0<p<2$, and $\mathcal R_{p,q}$ for $p\geq2$, from below by transport costs on $Q$ without taking an infimum over isometries.

The measure $\mu$ places a significant amount of mass on $a_0$ so that we get a simple lower bound on the probability that a coupling $\gamma$ transports $a_0$ to $S(a_0)$. Namely, for any $\gamma \in \Pi(\mu,\nu)$, define $\tilde\gamma := (\operatorname{id},S^{-1})_\# \gamma$ and let $(X,Y)\sim\tilde\gamma$. Then,
\begin{align*}
    \tilde\gamma(\{(a_0,a_0)\}) &= \mathbb P (X=a_0,\ Y=a_0)\ge1-\mathbb P(X\ne a_0)-\mathbb P(Y\ne a_0) \ge \frac{1}{3}.
\end{align*}
The exact constants in the construction are convenient choices; their role is to give $a_0$ mass greater than $1/2$ and to separate the reference points from one another and from the cube. Next, we would also like to show a lower bound on the probability that $\tilde\gamma$ transports $a_k$ to $a_k$. To this end, we consider how the mass $\mu$ places on $a_k$ splits,
\[
    \mu(\{a_k\}) = \tilde\gamma(\{(a_k,a_k)\}) + \underbrace{\sum_{\substack{0\leq j\leq d\\j \neq k}} \tilde\gamma(\{(a_k,a_j)\}) + \tilde\gamma(\{a_k\} \times Q)}_{=:\, L_k(\tilde\gamma)},
\]
and we bound $L_k(\tilde\gamma)$ from above. Recall from~\eqref{eq:K} that, with Euclidean distances on both factors,
\[
  J_{p,q}(\tilde\gamma)=\iint\bigl|\,\|x-x'\|^q-\|y- y'\|^q\,\bigr|^p\,d\tilde\gamma(x,y)\,d\tilde\gamma(x',y'),
\]
which equals $J_{p,q}(\gamma)$ when $S$ is an isometry. Since $\tilde\gamma(\{(a_0,a_0)\}) \geq \frac{1}{3}$, a draw $(x',y')\sim\tilde\gamma$ equals $(a_0,a_0)$ with probability at least $\frac{1}{3}$, and in this case we may lower bound $J_{p,q}(\tilde\gamma)$ by
\[
    J_{p,q}(\tilde\gamma) \geq \tilde\gamma(\{(a_0,a_0)\}) \int \bigl|\,\|x-a_0\|^q-\|y- a_0\|^q\,\bigr|^p\, d\tilde\gamma(x,y).
\]
Now suppose that $\tilde\gamma$ couples $a_k$ with $a_j$ for $j \neq k$. Then $\|a_k-a_0\|=r_k$ and $\|a_j-a_0\|=r_j$, and
\[
    J_{p,q}(\tilde\gamma) \geq \tilde\gamma(\{(a_0,a_0)\}) \cdot  \tilde\gamma(\{(a_k,a_j)\}) |r_k^q - r_j^q|^p.
\]
Hence, $\tilde\gamma(\{(a_k,a_j)\}) \leq \frac{ J_{p,q}(\tilde\gamma)}{\tilde\gamma(\{(a_0,a_0)\}) \cdot |r_k^q - r_j^q|^p}$. Since we have chosen each $r_k$ such that $r_k \neq r_j$ for $j \neq k$, the denominator is nonzero.

Similarly, suppose that $\tilde\gamma$ couples $a_k$ with $y \in Q$. Since $Q = \prod_{i=1}^d \left[ \frac{1}{8\sqrt{d}}, \frac{3}{16\sqrt{d}}\right]$, we have that $\|y\|\le 3/16<1/2<r_k$, and hence
\[
  J_{p,q}(\tilde\gamma) \geq \tilde\gamma(\{(a_0,a_0)\})\int_{\{a_k\}\times Q}\bigl|\,r_k^q-\|y\|^q\,\bigr|^p\,d\tilde\gamma(x,y) \geq \tilde\gamma(\{(a_0,a_0)\}) \cdot \tilde\gamma(\{a_k\}\times Q) \cdot \left(r_k^q- \left(\frac{3}{16} \right)^q\right)^p.
\]
Therefore, we get an upper bound $\tilde\gamma(\{a_k\}\times Q)\le \frac{J_{p,q}(\tilde\gamma)}{\tilde\gamma(\{(a_0,a_0)\}) \cdot \bigl(r_k^q-(3/16)^q\bigr)^{p}}$. Combining these upper bounds on $\tilde\gamma(\{(a_k,a_j)\})$ and $\tilde\gamma(\{a_k\} \times Q)$ gives an upper bound on $L_k(\tilde\gamma)$ and a lower bound for $\tilde\gamma(\{(a_k,a_k)\})$.

\paragraph{Step 2: a lower bound for the plug-in error.}
Given $X_1,\dots,X_n \iid \mu$ and $Y_1,\dots,Y_m \iid \nu$, let $\hat\mu_n = \frac{1}{n} \sum_{i=1}^n \delta_{X_i}$ and $\hat\nu_m = \frac{1}{m} \sum_{j=1}^m \delta_{Y_j}$ be the corresponding empirical measures. Set $\hat\nu_m':=(S^{-1})_\#\hat\nu_m$. If $0 < p < 2$, then $S(x) = x$ and hence $D_{p,q}^p(\mu,\nu) = D_{p,q}^p(\mu,\mu) = 0$. Thus $\E\left|D_{p,q}^p(\hat\mu_n,\hat\nu_m)-D_{p,q}^p(\mu,\nu)\right| =  \E |D_{p,q}^p(\hat\mu_n,\hat\nu_m')|$.
Now consider $p\geq2$, and set $\vartheta:=2^{-q}$. For $x,x',y,y'\in\R^d$, write
$u:=\|x-x'\|^q$ and $v:=\|y-y'\|^q$. Since $S$ scales distances by $1/2$, pulling the second
marginal back by $S^{-1}$ gives
\[
 D_{p,q}^p(\hat\mu_n,\hat\nu_m)
 =\inf_{\gamma\in\Pi(\hat\mu_n,\hat\nu_m')}
 \iint|u-\vartheta v|^p\,d\gamma(x,y)\,d\gamma(x',y').
\]
We use the following inequality, proved in Section~\ref{sec:proof-scalar}.

\begin{lemma}\label{lem:scalar}
There exists a constant $c_{p,q}>0$, possibly dependent on $p$ and $q$, such that, for all
$u,v\geq0$,
\begin{equation}\label{eq:scalar}
 |u-\vartheta v|^p \geq (1-\vartheta)^{p-1}(u^p-\vartheta v^p) +c_{p,q}(u+v)^{p-2}(u-v)^2.
\end{equation}
\end{lemma}

For a measure $\rho$, define
\[
 I_{p,q}(\rho):=\iint\|x-x'\|^{pq}\,d\rho(x)\,d\rho(x'),
\]
and, for $p\ge2$, set
$$
\widetilde J_{p,q}(\gamma):=\iint (u+v)^{p-2}(u-v)^2\,d\gamma(x,y)\,d\gamma(x',y').
$$
For probability measures $\rho_0,\rho_1$ on $\R^d$, set
$$
 \mathcal R_{p,q}(\rho_0,\rho_1):=\inf_{\gamma\in\Pi(\rho_0,\rho_1)}\widetilde J_{p,q}(\gamma).
$$
Integrating the scalar inequality and taking the infimum gives
\[
 D_{p,q}^p(\hat\mu_n,\hat\nu_m)
 \geq(1-\vartheta)^{p-1}\bigl[I_{p,q}(\hat\mu_n)-\vartheta I_{p,q}(\hat\nu_m')\bigr]
 +c_{p,q}\mathcal R_{p,q}(\hat\mu_n,\hat\nu_m').
\]
For the population pair $\mu, \nu$, the same inequality gives the lower bound
$(1-\vartheta)^p I_{p,q}(\mu)$, and the diagonal coupling of $\mu$ with itself attains it. Thus $D_{p,q}^p(\mu,\nu)=(1-\vartheta)^p I_{p,q}(\mu)$. Both $\hat\mu_n$ and $\hat\nu_m'$ are empirical measures from $\mu$, so that the diagonal terms vanish in the double sum defining $I_{p,q}(\hat\mu_n)$, yielding $\E I_{p,q}(\hat\mu_n)=\left(1-\frac1n\right)I_{p,q}(\mu)$ and $\E I_{p,q}(\hat\nu_m')=\left(1-\frac1m\right)I_{p,q}(\mu)$.
Subtracting $D_{p,q}^p(\mu,\nu)$, taking expectations, and using $|a-b|\ge a-b$ gives
\[
 \E| D_{p,q}^p(\hat\mu_n,\hat\nu_m) - D_{p,q}^p(\mu,\nu)|\geq (1-\vartheta)^{p-1}\left(\frac{\vartheta}{m}-\frac1n\right)I_{p,q}(\mu)
 +c_{p,q}\E\mathcal R_{p,q}(\hat\mu_n,\hat\nu_m'),
\]
which gives $\E| D_{p,q}^p(\hat\mu_n,\hat\nu_m) - D_{p,q}^p(\mu,\nu)|
 \geq c_{p,q}\E\mathcal R_{p,q}(\hat\mu_n,\hat\nu_m')-\frac{2^{pq}}{N}$ since $I_{p,q}(\mu)\leq2^{pq}$.
Thus
\[
 \E|D_{p,q}^p(\hat\mu_n,\hat\nu_m) - D_{p,q}^p(\mu,\nu)|\geq\begin{cases}
 \E D_{p,q}^p(\hat\mu_n,\hat\nu_m'),&0<p<2,\\
 c_{p,q}\E\mathcal R_{p,q}(\hat\mu_n,\hat\nu_m')-2^{pq}/N,&p\geq2.
 \end{cases}
\]

\paragraph{Step 3: comparison with transport costs on $Q$.}
For $s>0$ and probability measures $P_0,P_1$ supported on $Q$, define the transport cost $W_s^s(P_0,P_1) := \inf_{\pi \in \Pi(P_0,P_1)}\int \|x-y\|^s\,d\pi(x,y).$
The following lemma gives the comparison.

\begin{lemma}\label{lem:fixed-points}
Define the event
\begin{align*}
    \mathcal E_{n,m} = &\left\{\hat\mu_n(\{a_0\}),\, \hat\nu_m(\{S(a_0)\}) \geq \frac{5}{8} \right\} \cap  \bigcap_{k=1}^d \left\{\hat\mu_n(\{a_k\}),\, \hat\nu_m(\{S(a_k)\}) \geq \frac{1}{8d} \right\} \cap \\
    &\qquad \left\{\hat\mu_n(Q),\, \hat\nu_m(S(Q)) \geq \frac{1}{12} \right\}.
\end{align*}
On $\mathcal E_{n,m}$, set $P_{\mu,n}:=\frac{\hat\mu_n|_Q}{\hat\mu_n(Q)}$ and $P_{\nu,m}:=\frac{\hat\nu_m'|_Q}{\hat\nu_m'(Q)}$, the probability measures obtained by restricting to $Q$ and dividing by the mass there. There exists a constant $c=c(p,q,d)>0$ such that, on $\mathcal E_{n,m}$:
\begin{enumerate}[label=\textup{(\alph*)}]
\item For $0<p<2$, $D_{p,q}^p(\hat\mu_n,\hat\nu_m')\geq c\,W_p^p(P_{\mu,n},P_{\nu,m})$.
\item If $p\geq2$, then $\mathcal R_{p,q}(\hat\mu_n,\hat\nu_m')\geq c\,W_2^2(P_{\mu,n},P_{\nu,m})$.
\end{enumerate}
\end{lemma}

On $\mathcal E_{n,m}$, set $k:=n\hat\mu_n(Q)$ and $l:=m\hat\nu_m(S(Q))$.
If $i_1<\cdots<i_k$ are the indices for which $X_i\in Q$, and $j_1<\cdots<j_l$ are those for which $Y_j\in S(Q)$, then $P_{\mu,n}=\frac1k\sum_{r=1}^k\delta_{X_{i_r}}$ and similarly $P_{\nu,m}=\frac1l\sum_{r=1}^l\delta_{S^{-1}(Y_{j_r})}$. In particular, $n/12\leq k\leq n$ and $m/12\leq l\leq m$. We now apply Lemma~\ref{lem:fixed-points}. When $0<p<2$, the map $S$ preserves pairwise distances, so $D_{p,q}^p(\hat\mu_n,\hat\nu_m) \geq cW_p^p(P_{\mu,n},P_{\nu,m})$ on $\mathcal E_{n,m}$.
When $p\geq2$, Lemma~\ref{lem:fixed-points} gives $\mathcal R_{p,q}\bigl(\hat\mu_n,\hat\nu_m'\bigr) \geq cW_2^2(P_{\mu,n},P_{\nu,m})$ on $\mathcal E_{n,m}$. Here $\hat\nu_m'=m^{-1}\sum_{j=1}^m\delta_{2Y_j}$.

\paragraph{Step 4: a lower bound for the expected transport cost.}

We bound the expected transport cost on $Q$ from below.

\begin{lemma}\label{lem:finite-support}
Let $s>0$ and $N:=n\wedge m$. For the measures
$P_{\mu,n},P_{\nu,m}$ and the event $\mathcal E_{n,m}$
defined above, there exists $c_{s,d}>0$ such that
\begin{equation}\label{eq:conditional-transport}
 \E\!\left[
   \mathbf1_{\mathcal E_{n,m}}
   W_s^s(P_{\mu,n},P_{\nu,m})
 \right]
 \geq c_{s,d}N^{-s/d}\mathbb P(\mathcal E_{n,m}).
\end{equation}
The integrand is defined to be zero outside $\mathcal E_{n,m}$.
\end{lemma}

The full proof is given in Section~\ref{sec:proof-finite-support}; we give the idea here.
Suppose $n\leq m$, so $N=n$, and condition on $X_1,\dots,X_n$ and on which block $H_i$ contains
each pulled-back observation $S^{-1}(Y_j)$. This determines $\mathcal E_{n,m}$ and the indices of the
second-sample points in $Q$. On $\mathcal E_{n,m}$, the first-sample points in $Q$ form a nonempty set
$F$ of at most $n$ points, while the second-sample points in $Q$ are conditionally uniform there.
Let $v_d$ be the volume of the unit ball and set $r:=\ell(2v_dn)^{-1/d}$. The balls of radius $r$
centered at $F$ cover at most half of $Q$, so a uniform point in $Q$ has expected distance to $F$
raised to the power $s$ at least $r^s/2$. Any coupling to a measure supported on $F$ must pay at
least this distance cost. Taking conditional expectations and then expectations gives the lemma.

\paragraph{Step 5: a lower bound for $\mathbb P(\mathcal E_{n,m})$.}

We now prove that the conditions for Lemma~\ref{lem:fixed-points} hold with probability $\geq 1/2$ for sufficiently large $N$. To do this, we bound the probability that any empirical mass is smaller than required. Let $B\in\{\{a_0\},\dots,\{a_d\},Q\}$, set $b:=\mu(B)$, and let $\hat b:=\hat\mu_n(B)=n^{-1}\sum_{i=1}^n\mathbf1_{\{X_i\in B\}}$. Note that $\mathbf1_{\{X_i\in B\}}\sim\operatorname{Ber}(b)$, so $n\hat b\sim\operatorname{Bin}(n,b)$ and $\Var(\hat b) \leq \frac1{4n}$.

Set $t_B:=5/8$ if $B=\{a_0\}$, $t_B:=1/(8d)$ if $B=\{a_j\}$ for $1\leq j\leq d$, and $t_B:=1/12$ if $B=Q$.
In each case we have $b-t_B \geq 1/(24d)$, so Chebyshev's inequality gives
\[
 \mathbb P(\hat b<t_B)
 \leq\mathbb P(|\hat b-b|\geq b-t_B)
 \leq\frac{144d^2}{n}.
\]
Similarly, $\mathbb P(\hat\nu_m(S(B))<t_B)\leq144d^2/m$.
A union bound over the $d+2$ blocks and the two samples gives
$$
\mathbb P(\mathcal E_{n,m}^c)
\leq144d^2(d+2)(n^{-1}+m^{-1})
\leq\frac{288d^2(d+2)}{N},
$$
so $\mathbb P(\mathcal E_{n,m})\geq1/2$ for sufficiently large $N$.

\begin{proof}[Proof of Proposition~\ref{prop:lower-large-dim}]
If $0<p<2$, the above steps give
\begin{align*}
 \E\left|D_{p,q}^p(\hat\mu_n,\hat\nu_m)-D_{p,q}^p(\mu,\nu)\right|
 &=\E |D_{p,q}^p(\hat\mu_n,\hat\nu_m)|\\
 &\geq c\,\E\bigl[\mathbf1_{\mathcal E_{n,m}}W_p^p(P_{\mu,n},P_{\nu,m})\bigr]\\
 &\geq c' N^{-p/d}\cdot \mathbb P(\mathcal E_{n,m})
 \geq\frac{c'}2N^{-p/d}.
\end{align*}

If $p\geq2$, we have $d\geq2(p\wedge2)=4$. For sufficiently large $N$, the above steps give
\begin{align*}
 \E \left|D_{p,q}^p(\hat\mu_n,\hat\nu_m) - D_{p,q}^p(\mu,\nu) \right| &\geq c_{p,q} \E |\mathcal R_{p,q}(\hat\mu_n,\hat\nu_m')| - \frac{2^{pq}}{N} \\
 &\geq c\,\E\bigl[\mathbf1_{\mathcal E_{n,m}}W_2^2(P_{\mu,n},P_{\nu,m})\bigr] - \frac{2^{pq}}{N}\\
 &\geq c' N^{-2/d} \cdot \mathbb P(\mathcal E_{n,m}) - \frac{2^{pq}}{N} \geq\frac{c'}2N^{-2/d} - \frac{2^{pq}}{N} \\
 &\geq \frac{c'}4 N^{-2/d}.
\end{align*}
Indeed, $N^{-1}\leq N^{-1/2}N^{-2/d}$, so the last inequality holds whenever
$N\geq(4\cdot2^{pq}/c')^2$. Combining the two cases completes the proof.
\end{proof}

\subsection{Proof of Lemma~\ref{lem:scalar}}\label{sec:proof-scalar}
\begin{proof}
Recall that $\vartheta=2^{-q}\in(0,1)$. For $p=2$, the claim follows from the identity
\[
 (u-\vartheta v)^2=(1-\vartheta)(u^2-\vartheta v^2)+\vartheta(u-v)^2.
\]
For $p>2$, apply convexity of $t\mapsto|t|^p$ to
$u=(1-\vartheta)(u-\vartheta v)/(1-\vartheta)+\vartheta v$ to obtain
\begin{equation}\label{eq:scalar-convexity}
 |u-\vartheta v|^p
 \geq(1-\vartheta)^{p-1}(u^p-\vartheta v^p).
\end{equation}
By strict convexity, equality holds if and only if
$(u-\vartheta v)/(1-\vartheta)=v$, equivalently $u=v$. For $t\in[0,1]$, define
\[
 F(t):=|(1+\vartheta)t-\vartheta|^p
 -(1-\vartheta)^{p-1}\{t^p-\vartheta(1-t)^p\}.
\]
Applying \eqref{eq:scalar-convexity} with $u=t$ and $v=1-t$ gives $F(t)\geq0$, with equality if and only if $t=1/2$.
Direct differentiation gives $F'(1/2)=0$ and
$$
F''(1/2)=p(p-1)2^{4-p}\vartheta(1-\vartheta)^{p-2}>0.
$$
Taylor's formula therefore yields
$$
\lim_{t\to1/2}\frac{F(t)}{(2t-1)^2}
=\frac{F''(1/2)}8
=\frac{p(p-1)\vartheta(1-\vartheta)^{p-2}}{2^{p-1}}>0.
$$
Also, $F(t)/(2t-1)^2>0$ for $t\neq1/2$.
Thus, defining its value at $1/2$ by this limit gives a continuous, strictly positive function on $[0,1]$.
Its minimum $c_{p,q}$ is therefore positive, and
\[
 F(t)\geq c_{p,q}(2t-1)^2,\qquad t\in[0,1].
\]
Here the constant depends only on $p,q$, since $\vartheta=2^{-q}$.

For $u+v>0$, substitute $t=u/(u+v)$ and multiply by $(u+v)^p$ to obtain
\eqref{eq:scalar}. The case $u=v=0$ is immediate.
\end{proof}

\subsection{Proof of Lemma~\ref{lem:fixed-points}}\label{sec:proof-fixed-points}
\begin{proof}
Let $p>0$ and $\gamma \in \Pi(\hat\mu_n,\hat\nu_m')$. We separate the proof into the cases $0<p\leq2$ and $p>2$.

\medskip\noindent\emph{Case 1: $0<p<2 $.}
We now prove that $W_p^p(P_{\mu,n},P_{\nu,m})\leq CJ_{p,q}(\gamma)$ for every coupling in $\Pi(\hat\mu_n,\hat\nu_m')$. We first bound the mass that $\gamma$ assigns to $H_i\times H_j$ with $i\neq j$. On $\mathcal E_{n,m}$, both marginals assign mass at least $5/8$ to $a_0$, so
$
\gamma(\{(a_0,a_0)\})\geq5/8+5/8-1=1/4
$.
For each $i\in\{0,\dots,d+1\}$, we decompose the mass transported from $H_i$ into the mass that stays in $H_i$ and the mass transported to the other sets,
\[
    \hat\mu_n(H_i) = \gamma(H_i \times H_i) + \sum_{\substack{0\leq j\leq d+1\\j\ne i}} \gamma(H_i \times H_j) =:  \gamma(H_i \times H_i) + L_i(\gamma).
\]
An upper bound for $L_i(\gamma)$ gives a lower bound for $\gamma(H_i \times H_i)$. All sums over $i\ne j$ below use indices in $\{0,\ldots,d+1\}$. Define $L_\gamma := \sum_{i=0}^{d+1} L_i(\gamma)$ to be the total mass that $\gamma$ transports between distinct sets $H_i$ and $H_j$. By construction of the measure $\mu$, we have that
$1/8\le \|y\|\le 3/16$ for every $y\in Q$ and that
$\|a_k\|=r_k>1/2$ for $1\le k\le d$. For $(x,y)\in H\times H$, consider the following cases.

\begin{itemize}
\item If $x=a_0$ and $y\in H\setminus\{a_0\}$, then $y\in Q$
or $y=a_k$ for some $k=1,\ldots,d$, and $\|y\|\ge 1/8$. Thus,
$$
\bigl|\|x-x'\|^q-\|y-y'\|^q\bigr|
\ge \mathbf{1}_{\{(x',y')=(a_0,a_0)\}}\|y\|^q
\ge \left(\frac18\right)^q
\mathbf{1}_{\{(x',y')=(a_0,a_0)\}}.
$$

\item If $x=a_j$ and $y=a_k$ with $j\ne k$, then
$$
\bigl|\|x-x'\|^q-\|y-y'\|^q\bigr|
\ge \mathbf{1}_{\{(x',y')=(a_0,a_0)\}}
\bigl|\|a_j\|^q-\|a_k\|^q\bigr|
=|r_j^q-r_k^q|\,
\mathbf{1}_{\{(x',y')=(a_0,a_0)\}}.
$$

\item Consider $x=a_k$ for $k=1,\ldots,d$ and $y\in Q$.
Then $\|y\|\le 3/16<1/2<r_k$, and
$$
\bigl|\|x-x'\|^q-\|y-y'\|^q\bigr|
\ge \mathbf{1}_{\{(x',y')=(a_0,a_0)\}}
\bigl|\|a_k\|^q-\|y\|^q\bigr|
\ge \bigl|r_1^q-(3/16)^q\bigr|\,
\mathbf{1}_{\{(x',y')=(a_0,a_0)\}}.
$$
\end{itemize}
The same inequalities hold after interchanging $x$ and $y$.
Thus, for every $i\ne j$, $x\in H_i$, and $y\in H_j$, there is a constant $\delta > 0$ depending only on $q,d$ with $\bigl|\|x-x'\|^q-\|y-y'\|^q\bigr|
\ge \delta\,\mathbf{1}_{\{(x',y')=(a_0,a_0)\}},$ giving
\begin{equation}\label{eq:leakage-control}
\begin{aligned}
J_{p,q}(\gamma)
&\ge
\sum_{i\ne j}
\int_{H_i\times H_j}\int_{H\times H}
\mathbf{1}_{\{(x',y')=(a_0,a_0)\}}\delta^p
\,d\gamma(x',y')\,d\gamma(x,y)\\
&=
\gamma(\{(a_0,a_0)\})\delta^p
\sum_{i\ne j}\gamma(H_i\times H_j)\\
&\ge \frac{\delta^p}{4}L_\gamma.
\end{aligned}
\end{equation}

On $\mathcal E_{n,m}$, $\hat\mu_n(\{a_k\})\ge 1/(8d)$
for each $k=0,\ldots,d$, so
\begin{equation}\label{eq:anchor-mass}
\begin{aligned}
\frac{1}{8d}
&\le\hat\mu_n(\{a_k\})
=\gamma(\{(a_k,a_k)\})+L_k(\gamma)\\
&\le \gamma(\{(a_k,a_k)\})+L_\gamma.
\end{aligned}
\end{equation}

Now fix $x=(x_1,\ldots,x_d)\in Q$ and
$y=(y_1,\ldots,y_d)\in Q$. For $0\le k\le d$, a simple
computation shows that
$\|x-a_k\|,\|y-a_k\|\in[1/8,15/16]$.
By the mean value theorem applied to $t\mapsto t^{2/q}$,
$\bigl|\|x-a_k\|^2-\|y-a_k\|^2\bigr|\le C_q\bigl|\|x-a_k\|^q-\|y-a_k\|^q\bigr|$
for some sufficiently large $C_q>0$.
Since $a_k=r_ke_k$, expanding
$\|x-a_k\|^2$ and
$\|y-a_k\|^2$ gives
\begin{align*}
2r_k(x_k-y_k)
=\bigl[\|x-a_0\|^2-\|y-a_0\|^2\bigr]
-\bigl[\|x-a_k\|^2-\|y-a_k\|^2\bigr].
\end{align*}
Taking absolute values and using $2r_k\ge1$ yields
\begin{equation}\label{eq:anchor-distance-control}
\begin{aligned}
\|x-y\|
&\le \sum_{k=1}^d|x_k-y_k|
\le dC_q\sum_{k=0}^d
\bigl|\|x-a_k\|^q-\|y-a_k\|^q\bigr|.
\end{aligned}
\end{equation}
Finally, using the subadditivity of $t\mapsto t^p$ on $[0,\infty)$ when
$0<p\le1$ and H\"older's inequality when $p>1$ gives
\begin{equation}\label{eq:anchor-distance-control-p}
\|x-y\|^p
\le C_{p,q,d}\sum_{k=0}^d
\bigl|\|x-a_k\|^q-\|y-a_k\|^q\bigr|^p
\end{equation}
for a sufficiently large constant $C_{p,q,d}>0$. Restricting attention to the cube, we obtain
\begin{align*}
J_{p,q}(\gamma)
&\geq
\sum_{k=0}^d
\gamma(\{(a_k,a_k)\})
\int_{Q\times Q}
\bigl|\|x-a_k\|^q-\|y-a_k\|^q\bigr|^p
\,d\gamma(x,y).
\end{align*}
Also, $\|x-a_k\|,\|y-a_k\|\le15/16<1$ implies
$\bigl|\|x-a_k\|^q-\|y-a_k\|^q\bigr|^p\le1$.
By~\eqref{eq:anchor-mass},
\begin{align*}
&\frac{1}{8d}
\int_{Q\times Q}\sum_{k=0}^d
\bigl|\|x-a_k\|^q-\|y-a_k\|^q\bigr|^p
\,d\gamma(x,y)\\
&\qquad\le
\sum_{k=0}^d
\gamma(\{(a_k,a_k)\})
\int_{Q\times Q}
\bigl|\|x-a_k\|^q-\|y-a_k\|^q\bigr|^p
\,d\gamma(x,y)\\
&\qquad\quad+
L_\gamma
\int_{Q\times Q}\sum_{k=0}^d
\bigl|\|x-a_k\|^q-\|y-a_k\|^q\bigr|^p
\,d\gamma(x,y)\\
&\qquad\le J_{p,q}(\gamma)+(d+1)L_\gamma\\
&\qquad\le
\left(1+\frac{4(d+1)}{\delta^p}\right)J_{p,q}(\gamma).
\end{align*}
Combining this with \eqref{eq:anchor-distance-control-p}, we obtain
\begin{equation}\label{eq:cube-distortion-bound}
\int_{Q\times Q}\|x-y\|^p\,d\gamma(x,y)
\le
8dC_{p,q,d}
\left(1+\frac{4(d+1)}{\delta^p}\right)J_{p,q}(\gamma).
\end{equation}
This holds for every coupling $\gamma$ under consideration. The restricted measure $\gamma_Q:=\gamma|_{Q\times Q}$ may omit mass that leaves the cube $Q$ and gets paired with a delta mass under coupling $\gamma$. We construct a coupling of $P_{\mu,n}$ and $P_{\nu,m}$ from $\gamma_Q$.

Set $w:=\hat\mu_n(Q)$, $w':=\hat\nu_m(S(Q))$, $t:=\gamma(Q\times Q)$, and $W:=\max\{w,w'\}$. The amounts of mass sent between the cube and the reference points are
\[
 e:=w-t=\gamma(Q\times A),\qquad
 e':=w'-t=\gamma(A\times Q).
\]
They satisfy $e+e'\leq L_\gamma$ and $W-t=\max\{e,e'\}$.
Let $(\gamma_Q)_0,(\gamma_Q)_1$ be the first and second marginals of $\gamma_Q$. Note that $(\gamma_Q)_0\leq wP_{\mu,n}$ and $(\gamma_Q)_1\leq w'P_{\nu,m}$, so the measures
\[
 \kappa_0:=P_{\mu,n}-W^{-1}(\gamma_Q)_0,\qquad
 \kappa_1:=P_{\nu,m}-W^{-1}(\gamma_Q)_1
\]
are nonnegative. Both these measures have the same total mass, namely
\[
 r:=\kappa_0(Q)=\kappa_1(Q)=1-\frac tW
   =\frac{\max\{e,e'\}}W\leq\frac{L_\gamma}W\leq12L_\gamma.
\]
If $r>0$, define $\pi_Q :=W^{-1}\gamma_Q+r^{-1} \kappa_0 \otimes \kappa_1$; a simple computation shows that the first and second marginals of $\pi_Q$ are $P_{\mu,n}$ and $P_{\nu,m}$, respectively. If $r=0$, take $\pi_Q := W^{-1}\gamma_Q$. In either case, we have that $\pi_Q \in \Pi(P_{\mu,n},P_{\nu,m})$. The measure $r^{-1} \kappa_0 \otimes \kappa_1$ has total mass $r$ and is supported on $Q \times Q$ and also $W \geq 1/12$ on $\mathcal E_{n,m}$. Therefore,
\begin{equation}\label{eq:completion}
 \begin{aligned}
 W_p^p(P_{\mu,n},P_{\nu,m}) &\leq \int_{Q\times Q}\|x-y\|^p\,d\pi_Q(x,y) \\
 &\leq\frac1W\int_{Q\times Q}\|x-y\|^p\,d\gamma(x,y)+\diam(Q)^p r\\
 &\leq12\int_{Q\times Q}\|x-y\|^p\,d\gamma(x,y)+12\diam(Q)^pL_\gamma\\
 &\leq CJ_{p,q}(\gamma),
 \end{aligned}
\end{equation}
where the last inequality follows from \eqref{eq:leakage-control} and \eqref{eq:cube-distortion-bound}. The above argument applies to every $\gamma\in\Pi(\hat\mu_n,\hat\nu_m')$, so taking the infimum over all couplings gives
\[
     W_p^p(P_{\mu,n},P_{\nu,m}) \leq C\inf_{\gamma\in\Pi(\hat\mu_n,\hat\nu_m')}J_{p,q}(\gamma)=CD_{p,q}^p(\hat\mu_n,\hat\nu_m').
\]
This proves the first part of Lemma~\ref{lem:fixed-points}.

\medskip\noindent\emph{Case 2: $p\geq2$.} Fix $\gamma\in\Pi(\hat\mu_n,\hat\nu_m')$ and recall that
\[
 \widetilde J_{p,q}(\gamma)=\iint(u+v)^{p-2}(u-v)^2\,d\gamma(x,y)\,d\gamma(x',y'),
 \qquad u:=\|x-x'\|^q,\quad v:=\|y-y'\|^q.
\]
The integrand is zero when $u=v=0$.
We prove $W_2^2(P_{\mu,n},P_{\nu,m})\leq C\widetilde J_{p,q}(\gamma)$ for every such coupling.
Use the same sets $H_i$, number $\delta$, and mass $L_\gamma$ as in Case 1.
The mass at $(a_0,a_0)$ is again at least $1/4$.
Since $(u+v)^{p-2}(u-v)^2\geq|u-v|^p$, we have
$\widetilde J_{p,q}(\gamma)\geq J_{p,q}(\gamma)$, and~\eqref{eq:leakage-control} gives
\begin{equation}\label{eq:quadratic-cross-mass}
 L_\gamma\leq\frac{4\widetilde J_{p,q}(\gamma)}{\delta^p}.
\end{equation}
The bound~\eqref{eq:anchor-mass} also holds for every delta mass at $a_k$. By~\eqref{eq:anchor-distance-control},
\[
 \|x-y\|^2\leq K_2\sum_{k=0}^d
       \big|\|x-a_k\|^q-\|y-a_k\|^q\big|^2,
 \qquad K_2:=(dC_q)^2(d+1).
\]
We have that $\|x - a_k\| \geq 1/8$ and $\|y - a_k\| \geq 1/8$ for any two points $x,y \in Q$, so
\[
    \{\|x-a_k\|^q+\|y-a_k\|^q\}^{p-2} \big|\|x-a_k\|^q-\|y-a_k\|^q\big|^2 \geq b_2\big|\|x-a_k\|^q-\|y-a_k\|^q\big|^2,
\]
where $b_2:=(2\cdot8^{-q})^{p-2}>0$. Restricting $(x,y)$ to $Q\times Q$ and $(x',y')$ to $(a_i,a_i)$ for some $i = 0,1,\dots,d$ yields
\[
 \widetilde J_{p,q}(\gamma)
 \geq b_2\sum_{k=0}^d\gamma(\{(a_k,a_k)\})
       \int_{Q\times Q}\big|\|x-a_k\|^q-\|y-a_k\|^q\big|^2
       \,d\gamma(x,y).
\]
Use \eqref{eq:anchor-mass} as in the estimate preceding \eqref{eq:cube-distortion-bound}.
Since $|\|x-a_k\|^q-\|y-a_k\|^q|^2\leq1$,
\[
 \begin{aligned}
 \frac1{8d}\int_{Q\times Q}\sum_{k=0}^d
       \big|\|x-a_k\|^q-\|y-a_k\|^q\big|^2\,d\gamma(x,y)
 &\leq b_2^{-1}\widetilde J_{p,q}(\gamma)+(d+1)L_\gamma\\
 &\leq\left(b_2^{-1}+\frac{4(d+1)}{\delta^p}\right)\widetilde J_{p,q}(\gamma).
 \end{aligned}
\]
Consequently,
\begin{equation}\label{eq:quadratic-restricted-cost}
 \int_{Q\times Q}\|x-y\|^2\,d\gamma(x,y)
 \leq8dK_2\left(b_2^{-1}+\frac{4(d+1)}{\delta^p}\right)\widetilde J_{p,q}(\gamma).
\end{equation}
Again, this holds for every coupling $\gamma \in\Pi(\hat\mu_n,\hat\nu_m')$. Now take the same $\pi_Q$ constructed as before, but using $W = \max\{\hat\mu_n(Q),\hat\nu_m'(Q)\}$. Repeating the same argument as in Case 1 gives
\[
     W_2^2(P_{\mu,n},P_{\nu,m}) \leq12\int_{Q\times Q}\|x-y\|^2\,d\gamma(x,y)+12\diam(Q)^2L_\gamma \leq C\widetilde J_{p,q}(\gamma),
\]
by~\eqref{eq:quadratic-cross-mass} and~\eqref{eq:quadratic-restricted-cost}.
Taking the infimum over all couplings in $\Pi(\hat\mu_n,\hat\nu_m')$ gives
\[
     W_2^2(P_{\mu,n},P_{\nu,m}) \leq C\inf_{\gamma\in\Pi(\hat\mu_n,\hat\nu_m')}\widetilde J_{p,q}(\gamma)=C\mathcal R_{p,q}(\hat\mu_n,\hat\nu_m').
\]
This proves the second part of Lemma~\ref{lem:fixed-points} for $p\geq2$.
\end{proof}

\subsection{Proof of Lemma~\ref{lem:finite-support}}\label{sec:proof-finite-support}

\begin{proof}
Assume $n\leq m$; the case $m<n$ is symmetric. Write
$Z_j:=S^{-1}(Y_j)$ and condition on the $\sigma$-algebra $\mathcal H$ generated by
$X_1,\ldots,X_n$ and the indicators $\mathbf1_{\{Z_j\in H_i\}}$, $1\leq j\leq m$,
$0\leq i\leq d+1$. The event $\mathcal E_{n,m}$ and the set $\mathcal J:=\{j:Z_j\in Q\}$ are
$\mathcal H$-measurable. Conditionally on $\mathcal H$, the points $Z_j$ with $j\in\mathcal J$ are
independent and uniform on $Q$.

On $\mathcal E_{n,m}$, the set $F:=\{X_i:X_i\in Q\}$ is nonempty, $|F|\leq n$, and
$P_{\mu,n}$ is supported on $F$. Let $v_d$ be the volume of the unit ball in $\R^d$ and set
$r:=\ell(2v_dn)^{-1/d}$. The balls of radius $r$ centered at points of $F$ have total volume at
most $nv_dr^d=\ell^d/2$. Consequently,
\[
 \int_Q\min_{x\in F}\|z-x\|^s\,d\lambda_Q(z)
 \geq\frac{r^s}{2}=c_{s,d}n^{-s/d},
 \qquad c_{s,d}:=\frac{\ell^s}{2(2v_d)^{s/d}}.
\]
Every coupling of $P_{\mu,n}$ and $P_{\nu,m}$ has its first marginal supported on $F$, so
\[
 W_s^s(P_{\mu,n},P_{\nu,m})
 \geq\frac1{|\mathcal J|}\sum_{j\in\mathcal J}\min_{x\in F}\|Z_j-x\|^s
 \qquad\text{on }\mathcal E_{n,m}.
\]
Using the conditional uniform distribution and the fact that $\mathcal E_{n,m}\in\mathcal H$ gives
\[
 \E\!\left[\mathbf1_{\mathcal E_{n,m}}W_s^s(P_{\mu,n},P_{\nu,m})\mid\mathcal H\right]
 \geq\mathbf1_{\mathcal E_{n,m}}c_{s,d}n^{-s/d}.
\]
Taking expectations and using $n=N$ proves the claim.
\end{proof}

\endgroup

\section*{Acknowledgments}
The authors thank Promit Ghosal for providing valuable feedback on the exposition.

\bibliographystyle{alpha}
\bibliography{references}

@misc{leung2026sample,
  author = {Leung, Pui Kuen and Okada, Riku and Wong, Samuel Lok-Hei},
  title = {Sample Complexity for the {2-Gromov--Wasserstein} Distance},
  year = {2026},
  eprint = {2607.27514},
  archivePrefix = {arXiv},
  primaryClass = {math.ST},
  doi = {10.48550/arXiv.2607.27514},
  url = {https://arxiv.org/abs/2607.27514},
  note = {arXiv:2607.27514}
}

@misc{paliy2026empirical,
  author = {Paliy, Vasyl},
  title = {Empirical Convergence of Even-Order {Gromov--Wasserstein} Functionals},
  year = {2026},
  eprint = {2605.11108},
  archivePrefix = {arXiv},
  primaryClass = {math.PR},
  doi = {10.48550/arXiv.2605.11108},
  url = {https://arxiv.org/abs/2605.11108},
  note = {arXiv:2605.11108}
}

@book{santambrogio2015optimal,
  author = {Santambrogio, Filippo},
  title = {Optimal Transport for Applied Mathematicians:
    Calculus of Variations, {PDEs}, and Modeling},
  series = {Progress in Nonlinear Differential Equations and Their Applications},
  volume = {87},
  publisher = {Birkh{\"a}user},
  address = {Cham},
  year = {2015},
  doi = {10.1007/978-3-319-20828-2}
}

@article{guntuboyina2013covering,
  author = {Guntuboyina, Adityanand and Sen, Bodhisattva},
  title = {Covering Numbers for Convex Functions},
  journal = {IEEE Transactions on Information Theory},
  volume = {59},
  number = {4},
  pages = {1957--1965},
  year = {2013},
  eprint = {1204.0147},
  archivePrefix = {arXiv},
  url = {https://arxiv.org/abs/1204.0147}
}

@article{dudley1967sizes,
  author = {Dudley, Richard M.},
  title = {The Sizes of Compact Subsets of {Hilbert} Space and
    Continuity of {Gaussian} Processes},
  journal = {Journal of Functional Analysis},
  volume = {1},
  number = {3},
  pages = {290--330},
  year = {1967},
  doi = {10.1016/0022-1236(67)90017-1}
}

@article{zhang2024gromov,
  author = {Zhang, Z. and Goldfeld, Z. and Mroueh, Y. and Sriperumbudur, B. K.},
  title = {{Gromov--Wasserstein} Distances: Entropic Regularization,
    Duality and Sample Complexity},
  journal = {The Annals of Statistics},
  volume = {52}, number = {4}, pages = {1616--1645}, year = {2024},
  doi = {10.1214/24-AOS2406},
  url = {https://arxiv.org/abs/2212.12848}
}

@book{cannarsa2004semiconcave,
  author = {Cannarsa, P. and Sinestrari, C.},
  title = {Semiconcave Functions, {Hamilton--Jacobi} Equations, and Optimal Control},
  series = {Progress in Nonlinear Differential Equations and Their Applications},
  volume = {58}, publisher = {Birkh{\"a}user}, address = {Boston},
  year = {2004}, doi = {10.1007/b138356}
}

@book{vaart1996weak,
  author = {van der Vaart, A. W. and Wellner, J. A.},
  title = {Weak Convergence and Empirical Processes: With Applications to Statistics},
  series = {Springer Series in Statistics}, publisher = {Springer},
  address = {New York}, year = {1996}, doi = {10.1007/978-1-4757-2545-2}
}

@article{hundrieser2024empirical,
  author = {Hundrieser, S. and Staudt, T. and Munk, A.},
  title = {Empirical Optimal Transport Between Different Measures
    Adapts to Lower Complexity},
  journal = {Annales de l'Institut Henri Poincar{\'e}, Probabilit{\'e}s et Statistiques},
  volume = {60}, number = {2}, pages = {824--846}, year = {2024},
  doi = {10.1214/23-AIHP1369},
  url = {https://arxiv.org/abs/2202.10434}
}

@article{luxburg2004distance,
  author = {von Luxburg, U. and Bousquet, O.},
  title = {Distance-Based Classification with {Lipschitz} Functions},
  journal = {Journal of Machine Learning Research},
  volume = {5}, pages = {669--695}, year = {2004},
  url = {https://www.jmlr.org/papers/v5/luxburg04b.html}
}

@article{fournier2015rate,
  author = {Fournier, N. and Guillin, A.},
  title = {On the Rate of Convergence in {Wasserstein} Distance of the Empirical Measure},
  journal = {Probability Theory and Related Fields}, volume = {162},
  pages = {707--738}, year = {2015}, doi = {10.1007/s00440-014-0583-7},
  url = {https://arxiv.org/abs/1312.2128}
}

@article{berend2013sharp,
  author  = {Berend, Daniel and Kontorovich, Aryeh},
  title   = {A sharp estimate of the binomial mean absolute
             deviation with applications},
  journal = {Statistics \& Probability Letters},
  volume  = {83},
  number  = {4},
  pages   = {1254--1259},
  year    = {2013},
  doi     = {10.1016/j.spl.2013.01.023}
}

@article{bronshtein1976entropy,
  author = {Bronshtein, E. M.},
  title = {$\varepsilon$-Entropy of Convex Sets and Functions},
  journal = {Siberian Mathematical Journal},
  volume = {17},
  number = {3},
  pages = {393--398},
  year = {1976},
  doi = {10.1007/BF00967858}
}

@article{carl1999metric,
  author = {Carl, Bernd and Kyrezi, Ioanna and Pajor, Alain},
  title = {Metric Entropy of Convex Hulls in {Banach} Spaces},
  journal = {Journal of the London Mathematical Society},
  volume = {60},
  number = {3},
  pages = {871--896},
  year = {1999},
  doi = {10.1112/S0024610799008005}
}

@article{memoli2011gromov,
  author = {M{\'e}moli, Facundo},
  title = {{Gromov--Wasserstein} Distances and the Metric Approach to Object Matching},
  journal = {Foundations of Computational Mathematics},
  volume = {11},
  number = {4},
  pages = {417--487},
  year = {2011},
  doi = {10.1007/s10208-011-9093-5}
}

@book{petersen2016riemannian,
  author = {Petersen, Peter},
  title = {Riemannian Geometry},
  edition = {Third},
  series = {Graduate Texts in Mathematics},
  volume = {171},
  publisher = {Springer},
  address = {Cham},
  year = {2016},
  doi = {10.1007/978-3-319-26654-1}
}

@article{mantegazza2003hamilton,
  author = {Mantegazza, Carlo and Mennucci, Andrea Carlo},
  title = {{Hamilton--Jacobi} Equations and Distance Functions on {Riemannian} Manifolds},
  journal = {Applied Mathematics \& Optimization},
  volume = {47},
  number = {1},
  pages = {1--25},
  year = {2003},
  doi = {10.1007/s00245-002-0736-4},
  url = {https://cvgmt.sns.it/document/2286/}
}

@book{sturm2023space,
  author = {Sturm, Karl-Theodor},
  title = {The Space of Spaces: Curvature Bounds and Gradient Flows on the Space of Metric Measure Spaces},
  series = {Memoirs of the American Mathematical Society},
  volume = {290, no.~1443},
  publisher = {American Mathematical Society},
  address = {Providence, RI},
  year = {2023},
  doi = {10.1090/memo/1443},
  eprint = {1208.0434},
  archivePrefix = {arXiv}
}

@book{villani2009optimal,
  author = {Villani, C{\'e}dric},
  title = {Optimal Transport: Old and New},
  series = {Grundlehren der mathematischen Wissenschaften},
  volume = {338},
  publisher = {Springer},
  address = {Berlin},
  year = {2009},
  doi = {10.1007/978-3-540-71050-9}
}

@techreport{song2002upper,
  author      = {Song, Shuguang and Wellner, Jon A.},
  title       = {An Upper Bound for Uniform Entropy Numbers},
  institution = {Department of Statistics, University of Washington},
  type        = {Technical Report},
  number      = {406},
  year        = {2002},
  url         = {https://stat.uw.edu/research/tech-reports/upper-bound-uniform-entropy-numbers}
}

@article{dumont2025existence,
  author = {Dumont, T. and Lacombe, T. and Vialard, F.-X.},
  title = {On the Existence of {Monge} Maps for the
    {Gromov--Wasserstein} Problem},
  journal = {Foundations of Computational Mathematics},
  volume = {25}, pages = {463--510}, year = {2025},
  doi = {10.1007/s10208-024-09643-0},
  url = {https://arxiv.org/abs/2210.11945}
}

@misc{dandapanthula2025optimal,
  author = {Dandapanthula, Sanjit and Podkopaev, Aleksandr and
    Kasiviswanathan, Shiva Prasad and Ramdas, Aaditya and Goldfeld, Ziv},
  title = {Optimal Transportation and Alignment Between {Gaussian} Measures},
  year = {2025},
  eprint = {2512.03579},
  archivePrefix = {arXiv},
  primaryClass = {cs.LG},
  note = {arXiv:2512.03579},
  doi = {10.48550/arXiv.2512.03579},
  url = {https://arxiv.org/abs/2512.03579}
}

@article{duchi2019variance,
  author = {Duchi, John and Namkoong, Hongseok},
  title = {Variance-based Regularization with Convex Objectives},
  journal = {Journal of Machine Learning Research},
  year = {2019},
  volume = {20},
  number = {68},
  pages = {1--55},
  url = {https://www.jmlr.org/papers/v20/17-750.html}
}

@misc{rioux2026discrete,
  author = {Rioux, Gabriel and Marks, Joanna and Passeggeri, Riccardo and Goldfeld, Ziv},
  title = {Discrete {Gromov--Wasserstein} Duality: Algorithms and Isomorphism Testing},
  year = {2026},
  eprint = {2609.03094},
  archivePrefix = {arXiv},
  primaryClass = {math.ST},
  doi = {10.48550/arXiv.2609.03094},
  url = {https://arxiv.org/abs/2609.03094},
  note = {arXiv:2609.03094}
}

@misc{kato2025convergence,
  author = {Kato, Kengo and Wang, Boyu},
  title = {Convergence of empirical {Gromov--Wasserstein} distance},
  year = {2025},
  eprint = {2508.03985},
  archivePrefix = {arXiv},
  primaryClass = {math.ST},
  doi = {10.48550/arXiv.2508.03985},
  url = {https://arxiv.org/abs/2508.03985},
  note = {arXiv:2508.03985}
}

\newpage

\appendix

\section{Auxiliary results}\label{app:tools}

Throughout this appendix, unless stated otherwise, the distances are normalized as at the
beginning of Section~\ref{sec:reduction}, so that $\diam\X\le1$, $\diam\Y\le1$ and
$0\le K_{p,q}\le1$. Covers and covering numbers are as in Section~\ref{sec:notation}.

\subsection{Deterministic lemmas}\label{sec:det}

Let $L_K:=(p\vee1)q^{p\wedge1}$. We record the following estimates for
$K_{p,q}$ and $c_\gamma$, which we use to obtain Hölder bounds for the
potentials in Lemma~\ref{lem:pot} and for the corrected costs in
Section~\ref{sec:yside}.

\begin{lemma}[Hölder bounds for the kernel and the cost]\label{lem:reg}
For all $x,\tilde x,x'\in\X$ and $y,\tilde y,y'\in\Y$,
\[
  \bigl|K_{p,q}((x,y),(x',y'))-K_{p,q}((\tilde x,\tilde y),(x',y'))\bigr|\le L_K\bigl(\dX(x,\tilde x)^{p\wedge1}+\dY(y,\tilde y)^{p\wedge1}\bigr).
\]
Consequently, for every $\gamma\in\PP(\X\times\Y)$, the cost $c_\gamma$ takes values in $[0,2]$ and
\[
  |c_\gamma(x,y)-c_\gamma(\tilde x,\tilde y)|\le2L_K\bigl(\dX(x,\tilde x)^{p\wedge1}+\dY(y,\tilde y)^{p\wedge1}\bigr).
\]
\end{lemma}

\begin{proof}
We use two elementary inequalities. First, for $a,b\in[0,1]$,
\[
  |a^q-b^q|\le q|a-b|;
\]
this follows from the mean value theorem since $q\ge1$. Second, for $u,v\in[-1,1]$,
\[
  \bigl||u|^p-|v|^p\bigr|\le(p\vee1)|u-v|^{p\wedge1};
\]
for $p\ge1$ this holds because $|\cdot|^p$ is $p$-Lipschitz on $[-1,1]$, and for $p<1$ because
$\bigl||u|^p-|v|^p\bigr|\le||u|-|v||^p\le|u-v|^p$.

Let $u:=\dX(x,x')^q-\dY(y,y')^q$ and $\tilde u:=\dX(\tilde x,x')^q-\dY(y,y')^q$, both in $[-1,1]$. By
the first inequality and the triangle inequality,
$$|u-\tilde u|\le q\,|\dX(x,x')-\dX(\tilde x,x')|\le q\dX(x,\tilde x).$$
Hence, by the second inequality, replacing $x$ by $\tilde x$ changes $K_{p,q}$ by at most
$(p\vee1)\bigl(q\dX(x,\tilde x)\bigr)^{p\wedge1}=L_K\dX(x,\tilde x)^{p\wedge1}$.
Replacing $y$ by $\tilde y$ changes $K_{p,q}$ by at
most $L_K\dY(y,\tilde y)^{p\wedge1}$ in the same way, and the triangle inequality gives the first
claim. The claims for $c_\gamma$ follow by integrating against $\gamma$.
\end{proof}

\begin{lemma}[Hölder and sup-norm bounds for the potentials]\label{lem:pot}
Let $\gamma\in\PP(\X\times\Y)$ and $(f,g)\in\Phi_{\X,\Y}(c_\gamma)$. Then
$$
|f(x)-f(\tilde x)|
\le\min\{2,2L_K\dX(x,\tilde x)^{p\wedge1}\},
\qquad x,\tilde x\in\X,
$$
and $g$ satisfies the analogous bound on $\Y$. If $\inf_\X f=0$, then
$|f|\le2$ on $\X$ and $|g|\le4$ on $\Y$.
\end{lemma}

\begin{proof}
Let $x,\tilde x\in\X$, and let $y\in\Y$ attain the minimum
$f(\tilde x)=c_\gamma(\tilde x,y)-g(y)$. Since $0\le c_\gamma\le2$,
Lemma~\ref{lem:reg} gives
$$
f(x)-f(\tilde x)
\le c_\gamma(x,y)-c_\gamma(\tilde x,y)
\le\min\{2,2L_K\dX(x,\tilde x)^{p\wedge1}\}.
$$
Exchanging $x$ and $\tilde x$ gives the bound for $f$, and the bound
for $g$ follows in the same way. The last claim is obvious.
\end{proof}

\subsection{Dudley's entropy bound}\label{sec:prob}
We record the version of Dudley's entropy bound used in the proof of
Proposition~\ref{prop:potentials}.

\begin{lemma}[Dudley's entropy bound \citep{dudley1967sizes}]\label{lem:dudley}
Let $\mathcal G$ be a nonempty class of measurable functions on a measurable space $\mathcal Z$,
uniformly bounded by $B_0>0$. Let $\zeta_1,\dots,\zeta_k$ be i.i.d.\ with law $P$ and empirical measure $\hat P_k$,
and let $\bar N(\eps,\mathcal G):=\sup_QN(\eps,\mathcal G,L^2(Q))$, the supremum over finitely supported
$Q$. Then
\begin{equation}\label{eq:dudley}
  \E\sup_{e\in\mathcal G}\Bigl|\int e\,d(\hat P_k-P)\Bigr|\le C_{B_0}\inf_{\delta\in(0,1]}\Bigl\{k^{-1/2}+\delta+k^{-1/2}\int_{\delta/4}^{1}\sqrt{\log\bigl(2\bar N(\eps,\mathcal G)\bigr)}\,d\eps\Bigr\}.
\end{equation}
Here $\bar N(\eps,\mathcal G)$ may be replaced by $N(\eps,\mathcal G,\|\cdot\|_\infty)\ge\bar N(\eps,\mathcal G)$. If $d\ge1$, $a>0$ and
$\log\bar N(\eps,\mathcal G)\le C\eps^{-d/a}$ for $\eps\in(0,1]$, the right side is at most
$C'\psi_{d,a}(k)$ for $k\ge2$, where $C'$ depends only on $B_0$, $C$, $d$ and $a$.
\end{lemma}

\begin{proof}
The bound~\eqref{eq:dudley} follows from symmetrization \citep[Lemma~2.3.1]{vaart1996weak} and,
conditionally on the sample, the generalized entropy bound of von Luxburg and
Bousquet~\citep[Theorem~16]{luxburg2004distance} (see also \citep[Equation~(14)]{duchi2019variance})
applied to $\mathcal G\cup(-\mathcal G)$, whose covering numbers in $L^2(\hat P_k)$ are at most
$2\bar N(\eps,\mathcal G)$ and equal $1$ for $\eps\ge2B_0$. For the rate, take
$\delta=k^{-1/2}$ if $d\le2a$ and $\delta=k^{-a/d}$ if $d>2a$.
\end{proof}

\section{Proof of Lemma~\ref{lem:riemR}}\label{app:riem}

For a $C^2$ map $\psi$ from an open set $W\subseteq\R^d$ into $(M,g)$ and $z\in W$, let
$d\psi(z)$ denote the differential of $\psi$ at $z$, and write $\|d\psi(z)\|$ for its
operator norm from the Euclidean norm to $g$. Let
$\|\nabla d\psi(z)\|:=\sup_{\|e\|=1}\bigl\|D_s\partial_s\psi(z+se)|_{s=0}\bigr\|_g$ be the largest covariant
acceleration of the image of a unit-speed line through $z$. If $c(s):=\psi(z+sv)$ and $u$ is $C^2$ near
$c(s)$, then
\begin{equation}\label{eq:chain}
  (u\circ c)''(s)=\Hess u(\dot c,\dot c)+du(D_s\dot c),\qquad \|\dot c\|_g\le\|d\psi(z+sv)\|\,\|v\|,\quad \|D_s\dot c\|_g\le\|\nabla d\psi(z+sv)\|\,\|v\|^2 .
\end{equation}

\begin{proof}[Proof of Lemma~\ref{lem:riemR}]
\emph{Charts and Lipschitz bounds.}
Let $M_2:=\{x\in M:d_M(x,\X)\le2\}$, which is compact by
Hopf--Rinow. Let $\iota>0$ be a lower bound for the injectivity radius on $M_2$, and let $\kappa_M\ge0$
bound the sectional curvatures on $M_2$ from above. Fix $\rho\in(0,1/8]$ with $3\rho<\iota$ and
$3\rho\sqrt{\kappa_M}<\pi/2$. By compactness, choose finitely many charts $\psi_j:B_3\to M$, smooth
up to $\bar B_3$, with $\X\subseteq\bigcup_j\psi_j(B_1)$ and $\psi_j(B_3)\subseteq\{d_M(\cdot,\X)<1/2\}$.
Let $C_0\ge1$ and $C_2$ bound $\|d\psi_j\|$ and $\|\nabla d\psi_j\|$ on $\bar B_2$. Since $B_2$ is
convex, $d_M(\psi_jz,\psi_jz')\le C_0\|z-z'\|$. Moreover, $U\subseteq\{d_M(\cdot,\X)<1/2\}$ and $\diam\X\le1$
give $\diam(\X\cup U)\le2$. For the second-difference bounds, let $a\in\X$, $z\pm v\in B_2$,
$\|v\|\le\delta_0$, $c(s):=\psi_j(z+sv)$ and $w(s):=d_M(a,c(s))$ for $s\in[-1,1]$. By the chart Lipschitz bound,
$|w(s)-w(0)|\le C_0\|v\|\le\rho/8$.

\medskip\noindent\emph{Lower bound for $\Delta_v^2w_{a,j}$ near $a$.}
Let $w(0)<2\rho$, so that $c([-1,1])\subseteq B(a,3\rho)$. Since $3\rho<\iota$, the
function $r_a:=d_M(a,\cdot)$ is smooth on $B(a,3\rho)\setminus\{a\}$, with $\|dr_a\|_g=1$. Since moreover
$\sqrt{\kappa_M}\,r_a<\pi/2$ there, the Hessian comparison theorem for $\sec\le\kappa_M$
\citep[Chapter~6]{petersen2016riemannian} gives
\[
  \Hess r_a\ge\sqrt{\kappa_M}\cot(\sqrt{\kappa_M}\,r_a)\,(g-dr_a\otimes dr_a)\ge0\qquad\text{on }B(a,3\rho)\setminus\{a\},
\]
where $\sqrt{\kappa_M}\cot(\sqrt{\kappa_M}\,r)$ is read as $1/r$ when $\kappa_M=0$. By~\eqref{eq:chain},
$w''(s)\ge-C_2\|v\|^2$ whenever $c(s)\ne a$. We may assume $v\ne0$; then, since $\psi_j$ is injective,
$c(s)=a$ for at most one $s=s_0$. Let $F(s):=w(s)+C_2\|v\|^2s^2/2$. If $a\notin c((-1,1))$,
then $F$ is convex on $(-1,1)$ and, by continuity, on $[-1,1]$; this includes the case where
$c(s)=a$ at an endpoint. Otherwise $s_0\in(-1,1)$, and $F$ is convex on $[-1,s_0]$ and on
$[s_0,1]$. Since $w\ge0=w(s_0)$, the one-sided derivatives of $F$ satisfy
$F'_-(s_0)\le C_2\|v\|^2s_0\le F'_+(s_0)$, so $F$ is again convex on $[-1,1]$. Whether or
not the curve passes through $a$, the convexity of $F$ gives
\[
  \Delta_v^2w_{a,j}(z)=F(1)+F(-1)-2F(0)-C_2\|v\|^2\ge-C_2\|v\|^2 .
\]

\medskip\noindent\emph{Upper bound for $\Delta_v^2w_{a,j}$ away from $a$.}
For the upper bound when $w(0)>\rho$, the distance from $a$ need not be
smooth in a neighborhood of $c(0)$, since $c(0)$ may lie in the cut locus
of $a$ (see~\citep[Proposition~4.8(1)]{mantegazza2003hamilton} with
$K=\{a\}$). We instead construct a smooth upper support function using a
point on a minimizing geodesic from $a$. Set $x:=c(0)=\psi_j(z)$. Take a
unit-speed minimizing geodesic $\sigma$ from $a$ to $x$, and let
$a':=\sigma(w(0)-\rho/2)$, so that $d_M(a',x)=\rho/2$ and
$d_M(a,a')=w(0)-\rho/2$. Since $x$ lies within $1/2$ of $\X$ and
$\diam\X\le1$, we have $d_M(a',a)\le d_M(a,x)<3/2$, so $a'\in M_2$.
By the triangle inequality, $w(s)\le\zeta(s):=d_M(a,a')+d_M(a',c(s))$,
with equality at $s=0$, so
$\Delta_v^2w_{a,j}(z)\le\zeta(1)+\zeta(-1)-2\zeta(0)$. Along $c$ we have
$d_M(a',c(s))\in[3\rho/8,5\rho/8]$. The function
$(x',y)\mapsto d_M(x',y)$ is smooth on the open set
$\{(x',y):0<d_M(x',y)<\mathrm{inj}(x')\}$, which contains the compact set
$\{(x',y)\in M_2\times M:\rho/4\le d_M(x',y)\le3\rho/4\}$ because
$3\rho<\iota$. Thus the Hessian of $d_M(x',\cdot)$ is bounded on this
compact set by some $H_\rho$, while its gradient has norm $1$.
By~\eqref{eq:chain}, $\zeta''\le(H_\rho C_0^2+C_2)\|v\|^2$ on $[-1,1]$,
and therefore $\zeta(1)+\zeta(-1)-2\zeta(0)\le(H_\rho C_0^2+C_2)\|v\|^2$.
Hence, with $K_1:=H_\rho C_0^2+C_2$,
\[
  w_{a,j}(z)<2\rho\ \Longrightarrow\ \Delta_v^2w_{a,j}(z)\ge-K_1\|v\|^2,
  \qquad
  w_{a,j}(z)>\rho\ \Longrightarrow\ \Delta_v^2w_{a,j}(z)\le K_1\|v\|^2.
\]

\medskip\noindent\emph{Lipschitz and second-difference bounds for $t_{a,j}$.}
It remains to pass from $w_{a,j}$ to $t_{a,j}=w_{a,j}^q$. Let
\[
  C_t:=q2^{q-1}C_0,\qquad \lambda_q:=q(q-1)\max\{(7\rho/8)^{q-2},2^{q-2}\},\qquad
  K_t:=\lambda_qC_0^2+q2^{q-1}K_1 .
\]
Since $w_{a,j}$ is $C_0$-Lipschitz with values in $[0,2]$, and $r\mapsto r^q$ has derivative in
$[0,q2^{q-1}]$ on $[0,2]$, the function $t_{a,j}$ is $C_t$-Lipschitz. For the second differences,
write $w_\pm:=w_{a,j}(z\pm v)$ and $w_0:=w_{a,j}(z)$, so that $|w_\pm-w_0|\le\rho/8$. If $w_0<2\rho$,
Jensen's inequality, followed by the tangent-line bound
$r^q\ge w_0^q+qw_0^{q-1}(r-w_0)$ with $r=(w_++w_-)/2$, gives
\[
  \Delta_v^2t_{a,j}(z)
  \ge 2\Bigl(\tfrac{w_++w_-}{2}\Bigr)^q-2w_0^q
  \ge qw_0^{q-1}(w_++w_--2w_0)
  \ge -q2^{q-1}K_1\|v\|^2
  \ge -K_t\|v\|^2,
\]
where, when $q=1$ and $w_0=0$, $qw_0^{q-1}$ denotes the right derivative $1$ at zero. If $w_0>\rho$,
then $w_\pm,w_0\in[7\rho/8,2]$, where $r\mapsto r^q$ has derivative in $[0,q2^{q-1}]$ and second
derivative at most $\lambda_q$. Lemma~\ref{lem:comp} with $G(r)=r^q$ on this interval and
$\Delta_\sharp=\Delta_+$, the upper bound for $w_{a,j}$, and $|w_+-w_-|\le2C_0\|v\|$ give
$\Delta_v^2t_{a,j}(z)\le\lambda_qC_0^2\|v\|^2+q2^{q-1}K_1\|v\|^2=K_t\|v\|^2$.
\end{proof}

\end{document}